\documentclass[11pt]{article}

\usepackage{amsmath,amssymb,amsfonts,amsthm}
\usepackage{isomath}
\usepackage{graphicx}
\usepackage{tikz,pgfplots}
\usepackage{booktabs,multirow}
\usepackage{subcaption} 
\usepackage{hyperref}
\usepackage{algorithm2e}
\usepackage{float}
\usepackage{cancel}
\usepackage[autostyle]{csquotes} 
\usepackage{stackengine}
\usepackage{setspace}
\usepackage{placeins}
\usepackage{esint} 
\usepackage{stmaryrd} 
\usepackage{adjustbox}
\usepackage{pdflscape}
\usepackage{mathabx}
\usepackage{tikz-3dplot}
\usetikzlibrary{arrows.meta,plotmarks,patterns,external}
\pgfplotsset{compat=newest}
\pgfplotsset{every axis/.append style={
    scaled ticks = false, 
    tick label style={/pgf/number format/fixed}
}}

\newlength\figureheight 
\newlength\figurewidth 

\usepackage[a4paper, total={7.0in, 10in}]{geometry}

\newtheorem{remark}{Remark}

\newcommand*{\clbracket}{\{\mskip-5mu\{}
\newcommand*{\crbracket}{\}\mskip-5mu\}}
\newcommand*{\BBarchi}{\Bar{\Bar{\chi}}}

\title{\textbf{The high-order Hermite discrete correction function method for surface-driven electromagnetic problems}}

\author{
  Yann-Meing Law\thanks{Department of Mathematics and Industrial Engineering, Polytechnique Montr\'eal, C.P. 6079, succ. Centre-ville, Montr\'eal, QC H3C 3A7, Canada, Email: yann-meing.law@polymtl.ca}
}

\date{} 

\begin{document}

\maketitle

\begin{abstract}
The Hermite-Taylor method evolves all the variables and their derivatives through order $m$ in time to achieve a $2m+1$ order rate of convergence.
The data required at each node of the staggered Cartesian meshes used by this method makes the enforcement of boundary and interface conditions challenging.
In this work,
    we propose a novel correction function method, 
    referred to as the discrete correction function method,
    which provides all the data required 
    by the Hermite method
    near the surface where a condition is enforced.
The flexibility of the resulting Hermite-Taylor discrete correction function method is demonstrated by considering a 
    wide range of problems,
    including those with variable coefficients,
    discontinuous solutions at the interface,
    and generalized sheet transition conditions. 
Although the focus of this work is on Maxwell's equations,
    this high-order method can be adapted to other linear wave systems. 
Several numerical examples in two space dimensions are performed to verify the properties of the proposed method, 
    including long-time simulations.
\end{abstract}

\bigskip

\noindent\textbf{Keywords:} Hermite method,  Correction function method, High order,  Generalized sheet transition conditions, Interface conditions, Immersed boundary conditions

\medskip

\noindent\textbf{MSC Codes:} 35Q61, 65M70, 78A45


\section{Introduction}
The simulation of electromagnetic wave propagation in homogeneous media is well understood,
    and various high-order numerical methods have been developed to provide accurate and efficient simulations in this setting, 
    such as finite-difference \cite{Yee1966,Xie2002}, 
    discontinuous Galerkin \cite{Cockburn2001,Hesthaven2002},
    spectral \cite{Yang1997,Fan2002}, 
    and Hermite methods \cite{Goodrich2005},
    to name a few.
However, 
    in real-world applications, 
    physical domains involving different media and complex geometries often arise.
A typical example of an electronic device with such complexity is a waveguide composed of dielectric rods.

There are two main approaches to treating complex physical domains: 
    body-fitted meshes and immersed methods.
Body-fitted meshes use numerical methods 
    capable of handling unstructured meshes
    and curvilinear elements that conform to the geometry of the physical domain and interfaces between different media.
This approach is well-suited to finite element methods. 
In this work, 
    we consider the latter approach, 
    where interfaces and boundaries are immersed in a simple Cartesian mesh.
In this situation, 
    numerical methods must be modified locally to accurately enforce mathematical models on interfaces and boundaries to retain their order of accuracy, 
    thereby reducing the phase error for long-time simulations \cite{Kreiss1972}. 
Compatibility boundary conditions \cite{Henshaw1994,Henshaw2006}, 
    Fourier continuation \cite{Lyon2010}, 
    and  matched interface and boundary \cite{Zhao2004, Zhang2016} 
    methods are examples of immersed methods. 
In this work,
    we focus on enforcing various boundary and interface conditions for Hermite methods within this framework.
    
Hermite methods have shown remarkable properties for linear hyperbolic problems with periodic boundary conditions.
The Hermite-Taylor method introduced by Goodrich, 
    Hagstrom and Lorenz combines two staggered meshes, 
    a Hermite interpolation procedure in space and a Taylor method to evolve in time the variables and their derivatives through order $m$, 
    leading to a $2m+1$ order method \cite{Goodrich2005}. 
In $d$ space dimensions, 
    there are $(m+1)^d$ degrees of freedom at each node of the mesh. 
The interpolation procedure 
    and the time evolution of the data are performed locally and independently on each cell at each half time step,
    resulting in a high computation-to-communication ratio 
    and
    making this method suitable for GPU implementations \cite{Dye2014,Vargas2017,Appelo2018,Vargas2019}.
For linear symmetric hyperbolic problems, 
    the Hermite-Taylor method has a stability condition that depends 
    only on the maximum wave speed in the system and is independent of the order.
In \cite{Chen2012}, 
    this property is leveraged to develop a simple and yet efficient order-adaptive algorithm for problems characterized by highly localized features, 
    such as local pulses. 
Taking advantage of a leapfrog time-stepping method, 
    conservative Hermite methods have been developed for first-order linear wave systems \cite{Vargas2019},  
    the scalar wave equation \cite{Appelo2018} and dispersive Maxwell's equations \cite{AppeloHagstromLaw2025}. 
We refer the interested reader to \cite{Hagstrom2014} for a review of Hermite methods.

In the Hermite framework, 
    the enforcement of boundary conditions is challenging due to the data required on the boundary, 
    that is, 
    all the variables and their derivatives through order $m$. 
Three main approaches have been proposed to overcome this challenge. 
The first approach involves hybrid methods combining discontinuous Galerkin (DG) and Hermite methods for Maxwell's equations and the scalar wave equation \cite{Chen2014,Beznosov2021}. 
This approach leverages the flexibility of DG methods to handle body-fitted meshes near the boundary 
and couples them with the Hermite method, 
    which is applied on a Cartesian mesh away from the boundary. 
Due to the restrictive stability condition of  DG methods,
    local time-stepping procedures are needed to preserve the advantage of larger time steps in the regions where the Hermite method is applied. 
The second approach uses compatibility boundary conditions to
    provide the $(m+1)^d$ data on the boundary and considers curvilinear 
    meshes to handle complex geometries for the scalar wave equation \cite{Loya2025}. 
Note that boundary conditions were successfully enforced for both dissipative and conservative Hermite methods within this framework.
However,
    this approach is not directly applicable to first-order hyperbolic problems.
In this work, 
    we further improve upon the third approach, 
    namely the correction function method \cite{LawAppelo2025,LawAppeloHagstrom2025}.
    
The correction function method was first developed to provide high-order accurate finite-difference solutions to Poisson's equation with interface jump conditions \cite{Marques2011, Marques2017, Marques2019}. 
This method aims to provide local solutions,
    named correction functions, 
    of the unknown variables
    in the vicinity of the boundary and interfaces where the base method 
    used in the bulk,
    e.g., 
    a finite-difference scheme, cannot be applied. 
The system of equations describing the correction functions is composed of three 
    constraints that enforce
    the governing equations, 
    the boundary or interface conditions, 
    and the correction functions to match the base method's solution. 
A correction function is approximated by a polynomial function defined 
    in a small region, 
    called a local patch, 
    that encloses a part of the boundary or interface. 
Approximations of correction functions are obtained by 
    minimizing a functional that is a square measure of the residual of these constraints and involves volumetric and surface integrals. 
This approach has been extended to the scalar wave equation \cite{Abraham2018} and Maxwell's equations \cite{LawMarquesNave2020,LawNave2021,LawNave2022} to preserve the accuracy of finite-difference methods when immersed boundaries and interfaces are considered.
In \cite{Zhou2024},
    the authors propose a mesh-free collocation method to approximate the correction functions described by an elliptic problem without discontinuous coefficients.  
This leads to a solvable square linear system, 
    provided that the collocation points 
    are appropriately selected.

As mentioned before, 
    correction function methods have also been developed to provide the $(m+1)^d$ data required by the Hermite-Taylor method near  boundaries and interfaces when Maxwell's equations are considered. 
Unfortunately, 
    the matrices coming from the least squares approach to approximate the correction functions have large condition numbers that limit the order of accuracy of the Hermite-Taylor method. 
Moreover, 
    the computational cost of volumetric and surface numerical integrations increases with the order of the Hermite-Taylor correction function method.
The discrete correction function method proposed in this work 
    overcomes these two drawbacks. 

The paper is organized as follows.
In Section~\ref{sec:problem_definition}, 
    we introduce Maxwell's equations with different boundary conditions,
    Maxwell's interface problem,
    and simple generalized sheet transition conditions modeling the electromagnetic response of a flat metasurface.
The Hermite-Taylor method is briefly presented in Section~\ref{sec:hermite_taylor_method}.
Focusing on the two space dimensions simplification of Maxwell's equations using the transverse magnetic mode, 
    we then describe in detail the proposed discrete correction function method and its advantages in Section~\ref{sec:discrete_correction_function_method}.
Finally, 
    the properties of the Hermite-Taylor method combined with the discrete correction function method are verified using 
    several numerical examples, 
    including long-time simulations, 
    in Section~\ref{sec:num_examples}.


\section{Problem definition} \label{sec:problem_definition}

We are seeking numerical solutions to Maxwell's equations 
\begin{equation} \label{eq:Maxwell_equations}
		\mu \partial_t \mathbold{H} + \nabla\times \mathbold{E} = 0, 
        \qquad
		\epsilon \partial_t \mathbold{E} - \nabla\times \mathbold{H} = 0, \qquad
		\nabla\cdot(\epsilon\mathbold{E}) = 0, \qquad
		\nabla\cdot(\mu\mathbold{H}) = 0,
\end{equation}
	in a domain $\Omega \subset \mathbb{R}^2$ and a time interval $I=[t_0,t_f]$. 
Here, 
    $\mathbold{H}$ is the magnetic field, 
	$\mathbold{E}$ is the electric field, 
	$\mu$ is the magnetic permeability and $\epsilon$ is the electric permittivity.
The initial condition is given by $\mathbold{H}(\mathbold{x},t_0) = \mathbold{H}_0(\mathbold{x})$ and $\mathbold{E}(\mathbold{x},t_0) = \mathbold{E}_0(\mathbold{x})$ in $\Omega$.

Defining $\Gamma$ as the boundary of the domain $\Omega$ and $\mathbold{n}$ as the outward unit normal to $\Gamma$,    
    we consider the following boundary conditions on $\Gamma \times I$: 
    perfect electric conductor (PEC) 
\begin{equation} \label{eq:PEC_bnd_cdn}
    \mathbold{n} \times \mathbold{E} = 0, \qquad \mathbold{n}\cdot (\mu \mathbold{H}) = 0, 
\end{equation}
    perfect magnetic conductor (PMC) 
\begin{equation} \label{eq:PMC_bnd_cdn}
    \mathbold{n} \times \mathbold{H} = 0,  \qquad \mathbold{n}\cdot (\epsilon \mathbold{E}) = 0, 
\end{equation}
    and impedance boundary condition 
\begin{equation} \label{eq:impedance_bnd_cdn}
    \mathbold{E}\times\mathbold{n} + Z \mathbold{n}\times(\mathbold{H}\times\mathbold{n}) = 0.
\end{equation}
Here $Z=\sqrt{\mu/\epsilon}$ is the impedance. 
Further details and results on the well-posedness of Maxwell's equations with these boundary conditions can be found in \cite{Assous2018,Lindell2019}.

In this work, 
    we also address Maxwell's interface problem. 
We then consider that the domain $\Omega$ is subdivided into two subdomains $\Omega^+$ and $\Omega^-$, 
    and is such that $\Omega = \Omega^+ \cup \Omega^-$ and $\Omega^+\cap\Omega^- = \tilde{\Gamma}$, 
    where $\tilde{\Gamma}$ is a smooth interface between subdomains.
In this situation, 
    the physical parameters are considered to be piecewise constant functions, 
    that is 
\begin{equation} 
    \mu(\mathbold{x}) = \left\{ 
  \begin{array}{l}
    \mu^+ \quad \mbox{for} \quad \mathbold{x} \in \Omega^+, \\[1pt]
    \mu^- \quad \mbox{for} \quad \mathbold{x} \in \Omega^-,
  \end{array} \right. 
  \qquad \mbox{and} \qquad 
    \epsilon(\mathbold{x}) = \left\{ 
  \begin{array}{l}
    \epsilon^+ \quad \mbox{for} \quad \mathbold{x} \in \Omega^+, \\[1pt]
    \epsilon^- \quad \mbox{for} \quad \mathbold{x} \in \Omega^-.
  \end{array} \right. 
\end{equation}
The conventional interface conditions, 
    resulting from the juxtaposition of two different media,
    are given by 
\begin{equation} \label{eq:interface_cdns}
\tilde{\mathbold{n}}\times\llbracket \mathbold{H} \rrbracket = 0, \qquad   \tilde{\mathbold{n}}\cdot\llbracket \mu\,\mathbold{H} \rrbracket = 0, \qquad 
\tilde{\mathbold{n}}\times\llbracket \mathbold{E} \rrbracket = 0, \qquad  \tilde{\mathbold{n}}\cdot\llbracket \epsilon\,\mathbold{E} \rrbracket = 0, 
\end{equation}
    on $\tilde{\Gamma}\times I$.
Here $\llbracket \mathbold{A} \rrbracket = \mathbold{A}^+ - \mathbold{A}^-$ is the jump of the variable $\mathbold{A}$ on the interface, 
	$\mathbold{A}^\pm$ is the solution in $\Omega^\pm$ and 
	$\tilde{\mathbold{n}}$ is the unit normal to the interface $\tilde{\Gamma}$. 

Finally, 
    we also consider a metasurface, 
	modeled by a homogeneous zero-thickness sheet, 
	embedded in the domain $\Omega$ and located on the surface $\tilde{\Gamma}$. 
In this situation, 
    if we split $\mathbold{A}$ into an incident wave $\mathbold{A}^i$, 
	reflected wave $\mathbold{A}^r$ and transmitted wave $\mathbold{A}^t$, 
	we can write $\llbracket \mathbold{A} \rrbracket  = \mathbold{A}^{t} - (\mathbold{A}^{i} + \mathbold{A}^{r})$.
Note that subdomains $\Omega^{+}$ and $\Omega^-$ do not necessarily have different physical properties. 
The time-domain generalized sheet transition conditions (GSTC) modeling the electromagnetic response of a metasurface with only tangential polarizations require enforcing
\begin{equation} \label{eq:GSTC_time_domain} 
		\tilde{\mathbold{n}}\times\llbracket \mathbold{H} \rrbracket = \partial_t \mathbold{P}_s, 
        \qquad
        \tilde{\mathbold{n}} \cdot \llbracket \mu \mathbold{H} \rrbracket = - \mu_0 \nabla\cdot \mathbold{M}_s,
        \qquad 
		\tilde{\mathbold{n}}\times\llbracket \mathbold{E} \rrbracket = -\mu_0\partial_t \mathbold{M}_s, 
	\qquad 
        \tilde{\mathbold{n}} \cdot \llbracket \epsilon \mathbold{E} \rrbracket = - \nabla\cdot \mathbold{P}_s,
\end{equation}
	on $\tilde{\Gamma}\times I$.
Here $\mu_0$ and $\epsilon_0$ are the vacuum magnetic permeability and electric permittivity, 
	$\mathbold{P}$ is the surface electric polarization density, 
	$\mathbold{M}$ is the surface magnetic polarization density,
	$\nabla_s u = -\tilde{\mathbold{n}}\times(\tilde{\mathbold{n}}\times(\nabla u))$ is the surface gradient and 
	$\mathbold{A}_s = -\tilde{\mathbold{n}}\times(\tilde{\mathbold{n}}\times\mathbold{A})$ is the projection of the vector $\mathbold{A}$ on the surface $\tilde{\Gamma}$.
In this work, 
    we consider a monoanisotropic diagonal and uniform metasurface.
The surface polarization densities are then related to the electromagnetic fields by 
\begin{equation} \label{eq:polarization_general}
	\mathbold{P} = \epsilon_0 \BBarchi_{ee} \clbracket\mathbold{E}\crbracket, \qquad
		\mathbold{M} = \BBarchi_{mm} \clbracket\mathbold{H}\crbracket, \qquad \BBarchi_{ee} = \mbox{diag}(\chi_{ee}^{xx}, \chi_{ee}^{yy}, \chi_{ee}^{zz}), \qquad \BBarchi_{mm} = \mbox{diag}(\chi_{mm}^{xx}, \chi_{mm}^{yy}, \chi_{mm}^{zz}).
\end{equation}	
Here $\clbracket\mathbold{A}\crbracket = (\mathbold{A}^++\mathbold{A}^-)/2$ is the average of the variable $\mathbold{A}$ on the interface,
    and $\BBarchi_{ee}$ and $\BBarchi_{mm}$ are the effective surface susceptibilities. 
We refer the reader to \cite{Achouri2021} for further details on GSTC models for metasurfaces. 

\section{Hermite-Taylor method} \label{sec:hermite_taylor_method}

In this section, 
    we briefly introduce the Hermite-Taylor method in two space dimensions. 
Consider a domain $\Omega = [x_\ell,x_r]\times[y_b,y_t]$ and a time interval $I = [t_0, t_f]$. 
As mentioned before, 
    Hermite methods require two staggered meshes: the primal and the dual meshes. 
The nodes of the primal mesh are located at
\begin{equation}
    (x_i,y_j,t_n) = (x_\ell + i \Delta x,y_b + j \Delta y, t_0 + n\Delta t), 
\end{equation}
    for $i=0,\dots,N_x$, 
    $j=0,\dots,N_y$, 
    $n = 0,\dots,N_t$.
Here, 
    $N_{\circ}$ is the number of cells in the $\circ$ direction, 
    $\Delta x = (x_r-x_\ell)/N_x$, $\Delta y = (y_t-y_b)/N_y$ and $\Delta t = (t_f-t_0)/N_t$.
The nodes of the dual mesh correspond to the cell centers of the primal mesh, 
    that is,
\begin{equation}
    (x_{i+1/2}, y_{j+1/2}, t_{n+1/2}) = (x_i + \Delta x/2, y_j +\Delta y/2, t_n +\Delta t /2), 
\end{equation}
    for $i = 0,\dots, N_x-1$, 
    $j=0,\dots, N_y-1$,
    $n = 0,\dots, N_t-1$.

We now focus on the procedure to perform one time step.
The Hermite-Taylor method evolves the data from time $t_n$ to $t_{n+1}$ using two processes: a Hermite interpolation in space and a Taylor method in time. 
Assume periodic boundary conditions and that all the variables and their derivatives through order $m$ are known at the nodes of the primal mesh at time $t_n$.  
For each variable and for each cell $[x_i, x_{i+1}]\times[y_j,y_{j+1}]$ of the primal mesh,
    we compute the unique Hermite interpolant centered at $(x_{i+1/2},y_{j+1/2})$ of degree $(2m+1)^2$ satisfying $ \frac{\partial^{k+\ell} U}{\partial x^k \partial y^\ell}$ for $k,\ell = 0,\dots,m$, at each corner of the cell. 
Here, 
    $U$ is a given variable.
The Hermite polynomial coefficients are scaled derivatives of the variable $U$ located at the cell center at time $t_n$.

For each variable, 
    we now construct a space-time polynomial to evolve the variables and their derivatives 
    in time. 
Expanding each scaled coefficient of the interpolant in time results in the Hermite-Taylor polynomial 
\begin{equation}
P_{i+1/2,j+1/2}^U(x,y,t) = \sum_{k=0}^{2m+1} \sum_{\ell=0}^{2m+1} \sum_{s=0}^q c_{k,\ell,s}^U \bigg(\frac{x-x_{i+1/2}}{\Delta x}\bigg)^k \bigg(\frac{y-y_{j+1/2}}{\Delta y}\bigg)^\ell \bigg(\frac{t-t_n}{\Delta t}\bigg)^s.
\end{equation}
Here, 
    $q=2(2m+1)$ to have an exact Taylor expansion in time of the scaled coefficients.
From the Hermite interpolation step, 
    we know the coefficients $c_{k,\ell,0}$ for $k,\ell = 0, \dots, 2m+1$.
We then enforce the considered partial differential equations and their derivatives to compute the remaining coefficients $c_{k,\ell,s}$. 
This leads to a recursion relation on the scaled polynomial coefficients,
    e.g., 
    for the advection equation $\partial_t u - \partial_x u - \partial_y u = 0$,
    we obtain 
\begin{equation}
    \begin{aligned}
    & & c_{k,\ell,s}^u = \frac{(k+1)\Delta t}{s\Delta x} c_{k+1,\ell,s-1}^u + \frac{(\ell+1)\Delta t}{s\Delta y} c_{k,\ell+1,s-1}^u, \quad k,\ell = 0,\dots, 2m, \quad s=1,\dots, q,\\
    & &  c_{k,2m+1,s}^u = \frac{(k+1)\Delta t}{s\Delta x} c_{k+1,2m+1,s-1}^u, \quad k = 0,\dots, 2m, \quad s=1,\dots, q,\\
    & & c_{2m+1,\ell,s}^u = \frac{(\ell+1)\Delta t}{s\Delta y} c_{2m+1,\ell+1,s-1}^u, \quad \ell = 0,\dots, 2m, \quad s=1,\dots, q, \\
    & & c_{2m+1,2m+1,s}^u = 0, \quad s=1,\dots, q.
    \end{aligned}
\end{equation} 
Having now the values of  scaled coefficients $c_{k,\ell,s}^U$, 
    we evaluate the Hermite-Taylor polynomial and its derivatives through order $m$ at node $(x_{i+1/2},y_{j+1/2},t_{n+1/2})$ of the dual mesh.
This completes a half time step. 
The evolution of the data from the dual mesh at time $t_{n+1/2}$ to the primal mesh at time $t_{n+1}$ uses a similar procedure to achieve a full time step. 

For linear hyperbolic problems, 
    this method has a $2m+1$ order of accuracy, 
    and its stability condition depends only on the maximum wave speed in the problem and 
    is independent of its order. 
Further details and theoretical results on this method can be found in \cite{Goodrich2005}.
As mentioned before,
    enforcing boundary and interface conditions is challenging in the Hermite-Taylor method due to the required data and the use of Cartesian meshes. 
In the next section, 
    we propose a new correction function method to handle such conditions.

\section{Discrete correction function method} \label{sec:discrete_correction_function_method}

The correction function method (CFM) seeks approximations of the variables in a local space-time domain, 
    referred to as a local patch, 
    in the vicinity of the surface on which a condition is enforced. 
A local patch is required to enclose three components: a part of the surface, 
    nodes where the base method 
    cannot be directly applied, 
    referred to as CF nodes,
    and 
    nodes where the numerical solution coming from the base method is available,
    referred to as BM nodes. 
Here, 
    we assume the surface is not moving, 
    and therefore the local patches can be built in a pre-computation step.
    
In the CFM proposed in \cite{Marques2011}, 
    the functional to be minimized involves volumetric quadrature rules, 
    which are easier to derive for
	rectangular local patches. 
In this work, 
	the proposed discrete correction function method (DCFM) is free of any quadrature rules.
For given CF nodes, 
    a local patch is mainly used as a search area to identify
    the BM nodes and the closest part of the surface,
    as well as to scale the polynomial correction functions. 
Here, 
    we use square local patches for convenience.
That being said, 
	other geometries could be considered, 
    such as ellipsoidal regions \cite{Gabbard2023}.

\begin{figure}   
	\centering
	\begin{adjustbox}{max width=1.0\textwidth,center}
		\includegraphics[width=2.5in,trim={0cm 0cm 0cm 0cm},clip]{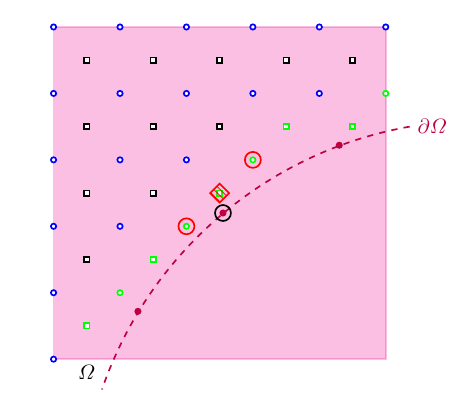} \hspace{1.0pt}
		\includegraphics[width=2.5in,trim={0cm 0cm 0cm 0cm},clip]{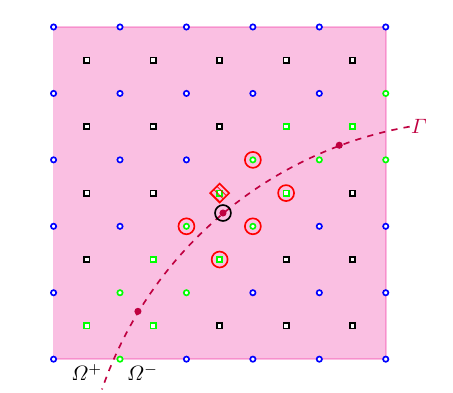}
	\end{adjustbox} 
       \caption{Examples of local patches. The left and right plots correspond to the situation with a boundary and an interface, 
       respectively. 
       The local patch is constructed for a node $u_j$ marked by a black circle.
       The green nodes represent the CF nodes,
       and those enclosed in red are associated with the local patch
       centered at the CF node marked by a dashed red square.
       The magenta region defines the spatial domain of the local patch.
       The primal and the dual BM nodes are illustrated by black squares and blue circles,
       respectively.}
\label{fig:examples_local_patches}
\end{figure}

The algorithm to generate local patches is detailed in \cite{LawAppeloHagstrom2025}. 
Here, 
    we briefly describe the procedure.
Assuming a parametrization of the surface, 
	we first discretize it in space with equidistant points denoted $u_j$.
The arc length between two consecutive points $u_j$ and $u_{j+1}$ is set to be the mesh size,
	denoted $h$, 
	in this work. 
We then seek the closest CF node to $u_j$ and center a square local patch at the spatial coordinate of this node, 
	denoted $\mathbold{x}_{c_j}$. 
In two space dimensions,
	the spatial domain of a local patch is then $\Omega_{c_j} = [x_{c_j}-\Delta x/2, x_{c_j}+\Delta x/2]\times[y_{c_j}-\Delta y/2, y_{c_j}+\Delta y/2]$ with $\Delta x = \Delta y = \beta h$ and $\beta$ is a given positive constant.
Here, 
    we choose $\beta = 6$.
Once the space domain of a local patch is defined, 
	we seek the BM nodes and the part of the surface within it. 
Finally, 
	for each CF node, 
	we assign the local patch for which its center $\mathbold{x}_{c_j}$ is the closest to the spatial coordinate of the CF node.
Typical local patches are illustrated in Fig.~\ref{fig:examples_local_patches}.
The time domain of a local patch depends on the mesh to which the enclosed CF nodes,
    where we seek a numerical solution at a given time, 
    belong. 
If the CF nodes belong to the primal mesh at time $t_{n+1}$, 
    the time domain is $I_{c_j} = [t_{n+1/2}, t_{n+1}]$, 
    while it is $I_{c_j} = [t_n, t_{n+1/2}]$ when 
    the CF nodes belong to the dual mesh at time $t_{n+1/2}$.
In both situations,
    $\Delta t = t_{n+1} - t_{n+1/2} = t_{n+1/2} - t_n$. 
    
Now that local patches are defined, 
	we compute the correction functions in such a way that three constraints are weakly enforced. 
The first constraint is to enforce the governing equations, 
    that is, 
    Maxwell's equations in this work.
The second constraint is to enforce the boundary, 
    interface or generalized sheet transition conditions on the surface.
The last constraint is to match the correction functions to the numerical solution obtained with the base method,
    here the Hermite-Taylor method.
    
The key idea of the novel correction function method is to directly constrain polynomial coefficients of the correction functions,
	hence the name discrete correction function method. 
As shown in this section, 
	this approach has two main advantages.
First, 
	we avoid the normal form of the least squares problem and therefore significantly reducing the condition number of matrices used in the minimization procedure. 
Second,
	this approach is free of any volumetric and surface quadrature rules, 
	and considers directly approximations of variables located at the BM nodes,
	making the DCFM more efficient than other correction function methods. 
In the following, 
	we first describe the method for the embedded boundary case. 
Afterwards, 
    we provide additional details for the spatially variable coefficients, 
    interface and metasurface cases.

Assume a local patch belonging to the boundary $\Gamma$. 
In two space dimensions, 
    we are seeking correction functions into the space of tensor-product polynomials of degree at most $d$, 
    denoted $\mathbb{Q}^d$, 
    in each independent variable,
    leading to polynomial correction functions centered at $(\mathbold{x}_c,t_c)$ of the form
\begin{equation}
    P^f_d(\mathbold{x},t) = \sum_{k=0}^d\sum_{\ell=0}^d\sum_{s=0}^d c_{k,\ell,s}^f \xi^k\eta^\ell\zeta^s, \qquad 
    \xi = \frac{x-x_c}{\Delta x}, \qquad \eta = \frac{y-y_c}{\Delta y}, \qquad \zeta = \frac{t-t_c}{\Delta t}.
\end{equation}
Here, 
    $f$ is a component of the electromagnetic fields $\mathbold{H}$ or $\mathbold{E}$,
 and $d>0$ to enforce constraints involving first-order derivatives, 
    such as in Maxwell's equations.
We choose $\mathbold{x}_c = (x_c,y_c)$ to be the center of the local patch and $t_c$ to be the time where we seek to evaluate the correction functions. 
The total number of coefficients to determine for the electromagnetic fields is $3(d+1)^3$.

\begin{remark}
It is important that all the equations that constrain the correction functions share the same dimensional unit to ensure a dimensionally consistent system of equations and to provide the scaling associated with different equations.
To do so, 
    we consider the reference quantities 
\begin{equation}
\mathbold{x} = L_0 \hat{\mathbold{x}}, \qquad t = \frac{L_0}{c_0} \hat{t}, \qquad \epsilon = \epsilon_0 \hat{\epsilon}, \qquad \mu = \mu_0 \hat{\mu}, \qquad \mathbold{H} = H_0 \hat{\mathbold{H}}, \qquad \mathbold{E} = Z_0 H_0 \hat{\mathbold{E}}. 
\end{equation}
Here, 
    $\hat{a}$ is the non-dimensional variable associated with $a$, 
    $L_0$ is the reference length, 
    $\mu_0$ is the magnetic permeability in free space, 
    $\epsilon_0$ is the electric permittivity in free space,
    $c_0 = 1/\sqrt{\epsilon_0\mu_0}$, 
    $Z_0 = \sqrt{\mu_0/\epsilon_0}$ and $H_0$ is the reference magnetic strength.
In the following, 
    we ensure that all the different terms in the system of equations for the correction function coefficients have the same dimensional unit as $Z_0 H_0$.
We also choose $L_0 = \min\{\Delta x, \Delta y\}$.
\end{remark}

\subsection{Governing equations constraints} \label{sec:gov_eq_constraints}

Let us consider a 2-D simplification of Maxwell's equations using the transverse magnetic (TM$_z$) mode with the appropriate scaling 
\begin{align}
    L_0(\mu\partial_t H_x + \partial_y E_z) =&\,\, 0, \label{eq:Faraday_law_1}\\
    L_0(\mu \partial_t H_y -\partial_x E_z) =&\,\, 0, \label{eq:Faraday_law_2}\\ Z_0L_0(\epsilon\partial_t E_z -\partial_x H_y + \partial_y H_x) =&\,\, 0, \label{eq:Ampere_Maxwell_law}\\
    c_0L_0(\mu \partial_x H_x + \mu \partial_y H_y) =&\,\, 0. \label{eq:divergence_free}
\end{align}
Substituting the electromagnetic fields with their polynomial correction function approximations in the first equation of Faraday's law \eqref{eq:Faraday_law_1} leads to 
\begin{equation}
    L_0(\mu\partial_t P_d^{H_x} + \partial_y P_d^{E_z}) = 0. 
\end{equation}
Expanding the polynomials and matching the coefficients gives
\begin{equation} \label{eq:cf_poly_step_1}
    \begin{aligned}
    & L_0\sum_{k=0}^d\sum_{\ell=0}^{d-1}\sum_{s=0}^{d-1} \bigg(\frac{\mu(s+1)}{\Delta t}c_{k,\ell,s+1}^{H_x} + \frac{(\ell+1)}{\Delta y}c_{k,\ell+1,s}^{E_z}\bigg)\xi^k\eta^\ell\zeta^s + L_0\sum_{k=0}^d\sum_{s=0}^{d-1} \frac{\mu(s+1)}{\Delta t}c_{k,d,s+1}^{H_x} \xi^k\eta^d\zeta^s  \\ 
    & + L_0\sum_{k=0}^d\sum_{\ell=0}^{d-1}\frac{(\ell+1)}{\Delta y}c_{k,\ell+1,d}^{E_z} \xi^k\eta^\ell\zeta^d = 0.  
    \end{aligned}
\end{equation}
For the polynomial in \eqref{eq:cf_poly_step_1} to vanish uniformly, 
    all its terms should be equal to zero, 
    and
    we obtain the following equations for the coefficients of $H_x$ and $E_z$
\begin{align}
   & \frac{\mu L_0}{\Delta t}(s+1)c_{k,\ell,s+1}^{H_x} + \frac{L_0}{\Delta y}(\ell+1)c_{k,\ell+1,s}^{E_z} = 0, \qquad k=0,\dots,d, \qquad \ell,s = 0,\dots,d-1, \label{eq:constraints_Faraday_1} \\
   & \frac{\mu L_0}{\Delta t}   (s+1) c_{k,d,s+1}^{H_x} = 0, \qquad k = 0,\dots,d, \qquad s = 0,\dots, d-1, \label{eq:constraints_Faraday_2}\\
   & \frac{L_0}{\Delta y} (\ell+1) c_{k,\ell+1,d}^{E_z} = 0, \qquad k = 0,\dots,d, \qquad \ell = 0,\dots, d-1. \label{eq:constraints_Faraday_3}
\end{align}
Repeating the same procedure for the second equation in Faraday's law \eqref{eq:Faraday_law_2},
    we obtain
\begin{align}
    & \frac{\mu L_0}{\Delta t}  (s+1) c_{k,\ell,s+1}^{H_y} - \frac{L_0}{\Delta x} (k+1) c_{k+1,\ell,s}^{E_z} = 0, \qquad k,s= 0, \dots, d-1, \qquad \ell = 0, \dots, d, \label{eq:constraints_Faraday_4} \\
    & \frac{\mu L_0}{\Delta t}  (s+1) c_{d,\ell,s+1}^{H_y} = 0, \qquad \ell = 0,\dots,d, \qquad s = 0, \dots, d-1, \label{eq:constraints_Faraday_5} \\
    &\frac{L_0}{\Delta x}(k+1) c_{k+1,\ell,d}^{E_z} = 0, \qquad k=0,\dots,d-1, \qquad \ell = 0,\dots,d. \label{eq:constraints_Faraday_6}
\end{align}
Considering now Amp\`ere-Maxwell law \eqref{eq:Ampere_Maxwell_law}, 
    we have the following equations for the coefficients 
\begin{align}
    & \frac{\epsilon Z_0L_0}{\Delta t} (s+1) c_{k,\ell,s+1}^{E_z} - \frac{Z_0L_0}{\Delta x} (k+1) c_{k+1,\ell,s}^{H_y} + \frac{Z_0L_0}{\Delta y} (\ell+1) c_{k,\ell+1,s}^{H_x} = 0, \qquad k,\ell,s= 0, \dots, d-1, \label{eq:constraints_Ampere_Maxwell_1}\\
    & \frac{\epsilon Z_0 L_0}{\Delta t} (s+1) c_{k,d,s+1}^{E_z} - \frac{Z_0L_0}{\Delta x} (k+1) c_{k+1,d,s}^{H_y} = 0, \qquad k,s = 0,\dots,d-1, \label{eq:constraints_Ampere_Maxwell_2}\\
    &\frac{\epsilon Z_0 L_0}{\Delta t} (s+1) c_{d,\ell,s+1}^{E_z} + \frac{Z_0 L_0}{\Delta y} (\ell+1) c_{d,\ell+1,s}^{H_x} = 0, \qquad \ell,s= 0, \dots, d-1, \label{eq:constraints_Ampere_Maxwell_3}\\
    & - \frac{Z_0 L_0}{\Delta x} (k+1) c_{k+1,\ell,d}^{H_y} + \frac{Z_0 L_0}{\Delta y} (\ell+1) c_{k,\ell+1,d}^{H_x} = 0, \qquad k,\ell= 0, \dots, d-1, \label{eq:constraints_Ampere_Maxwell_4}\\
    & \frac{\epsilon Z_0 L_0}{\Delta t} (s+1) c_{d,d,s+1}^{E_z} = 0, \qquad s= 0, \dots, d-1, \label{eq:constraints_Ampere_Maxwell_5}\\ 
    & - \frac{Z_0 L_0}{\Delta x} (k+1) c_{k+1,d,d}^{H_y} = 0, \qquad k= 0, \dots, d-1, \label{eq:constraints_Ampere_Maxwell_6}\\
    & \frac{Z_0 L_0}{\Delta y} (\ell+1) c_{d,\ell+1,d}^{H_x} = 0, \qquad \ell= 0, \dots, d-1. \label{eq:constraints_Ampere_Maxwell_7}
\end{align}
Finally, 
    the divergence-free constraint \eqref{eq:divergence_free} leads to 
\begin{align}
    & \frac{\mu c_0 L_0}{\Delta x} (k+1) c_{k+1,\ell,s}^{H_x} + \frac{\mu c_0 L_0}{\Delta y} (\ell+1) c_{k,\ell+1,s}^{H_y} = 0, \qquad k,\ell = 0, \dots, d-1, \qquad s = 0,\dots,d, \label{eq:constraints_div_1} \\
    & \frac{\mu c_0 L_0}{\Delta x} (k+1) c_{k+1,d,s}^{H_x} = 0, \qquad k =0,\dots,d-1, \qquad s =0,\dots,d, \label{eq:constraints_div_2} \\
    & \frac{\mu c_0 L_0}{\Delta y} (\ell+1) c_{d,\ell+1,s}^{H_y} = 0, \qquad \ell =0,\dots,d-1, \qquad s =0,\dots,d. \label{eq:constraints_div_last}
\end{align}
We then have $4d^3+12d^2+9d$ equations coming from the enforcement of the governing equations to satisfy.

\subsubsection{Strongly enforcing the divergence-free constraint}

In the discrete correction function method setting, 
    it is straightforward to enforce exactly the divergence-free constraint \eqref{eq:divergence_free} by choosing a divergence-free polynomial space.
Isolating the coefficients for the $x$-component of the magnetic field in equations \eqref{eq:constraints_div_1} leads to 
\begin{equation} \label{eq:simplification_div_free_cdn}
    c_{k,\ell,s}^{H_x} = -\frac{\Delta x (\ell+1)}{\Delta y k} c_{k-1,\ell+1,s}^{H_y}, \qquad k=1,\dots,d, \qquad \ell = 0, \dots, d-1, \qquad s = 0,\dots,d.
\end{equation}
From equations \eqref{eq:constraints_div_2} and \eqref{eq:constraints_div_last}, 
    we obtain 
\begin{align}
    & c_{k,d,s}^{H_x} = 0, \qquad k =1,\dots,d, \qquad s =0,\dots,d,  \\
    & c_{d,\ell,s}^{H_y} = 0, \qquad \ell =1,\dots,d, \qquad s =0,\dots,d.
\end{align}
Considering these equations directly in the definition of the polynomials 
    approximating the magnetic field results in seeking the associated correction functions within a divergence-free polynomial space,
    reducing the total number of coefficients for the electromagnetic fields to $2d^3+6d^2+7d+3$.
For large $d$, 
    the total number of coefficients is reduced by roughly one-third, 
    making this approach computationally more efficient.
Moreover,
    we also reduce the number of equations required to enforce Maxwell's equations.
In addition to removing equations \eqref{eq:constraints_div_1}-\eqref{eq:constraints_div_last} coming from the divergence-free constraint, 
    equations \eqref{eq:constraints_Faraday_2},
    \eqref{eq:constraints_Faraday_5},
    \eqref{eq:constraints_Ampere_Maxwell_6} and 
    \eqref{eq:constraints_Ampere_Maxwell_7} are simplified respectively to 
\begin{align}
& \frac{\mu L_0}{\Delta t}   (s+1) c_{0,d,s}^{H_x} = 0,  \qquad s = 0,\dots, d-1, \\ 
& \frac{\mu L_0}{\Delta t}  (s+1) c_{d,0,s+1}^{H_y} = 0, \qquad s = 0, \dots, d-1, \\
& - \frac{Z_0 L_0}{\Delta x} (k+1) c_{k+1,d,d}^{H_y} = 0, \qquad k= 0, \dots, d-2, \label{eq:constraints_Ampere_Maxwell_6_modify}\\ 
& \frac{Z_0 L_0}{\Delta y} (\ell+1) c_{d,\ell+1,d}^{H_x} = 0, \qquad \ell= 0, \dots, d-2. \label{eq:constraints_Ampere_Maxwell_7_modify}
\end{align}
We then have a total of $3d^3+7d^2+7d-2$ equations coming from Maxwell's equations to satisfy when divergence-free correction function polynomials are considered. 
In the remainder of this work, 
    we always consider divergence-free polynomial spaces,
    although we do not explicitly modify the notation of correction function coefficients for the sake of legibility.
    
\begin{remark}
Consider a  divergence constraint that is non-trivial, 
    e.g., $\nabla\cdot(\epsilon \mathbold{E}) = \rho(\mathbold{x})$ in TE$_z$ mode.
We could still reduce the total number of coefficients in the correction function polynomials and the number of equations for these coefficients if the scalar charge density $\rho(\mathbold{x})$ is accurately approximated locally by a polynomial. 
As an example, 
    if $\rho$ is a positive constant, 
    we obtain similar equations as in \eqref{eq:simplification_div_free_cdn}, 
    but for coefficients $c_{k,\ell,s}^{E_x}$ and $c_{k,\ell,s}^{E_y}$, 
    except for the subscript $(k,\ell,s) = (1,0,0)$, 
    for which we have 
\begin{equation}
    c_{1,0,0}^{E_x} = \rho - \frac{\Delta x}{\Delta y}c_{0,1,0}^{E_y}. 
\end{equation}
Unfortunately, 
    the divergence constraint is not enforced exactly anymore
    since the terms involving $\rho$ are part of the right-hand side of the overdetermined linear system of equations enforcing the governing equations.
In this situation,
    the divergence constraint is weakly enforced using a least squares approach,
    as described in subsection~\ref{sec:solving_linear_system_bnd_cdns}.
\end{remark}

\subsection{Matching the base method's solution constraints}

We now constrain the correction functions to match the numerical solution coming from the base method. 
Assume the numerical solution is available at the BM nodes within a local patch for a variable $f$, 
    that is 
    $f^*(\mathbold{x}^f_{i},t^f_{i})$ for $i=1,\dots, N^f$. 
As an example, 
    if we consider the CF nodes to be updated belong to the primal mesh at time $t_{n+1}$,
    the time domain of their local patches includes the Hermite-Taylor numerical solution at the dual nodes at time $t_{n+1/2}$ and at the primal nodes at $t_{n+1}$.

Having identified the BM nodes within a local patch, 
    we therefore enforce the following equations 
\begin{align}
	& \frac{\omega_b Z_0}{N^{H_x}} P_d^{H_x}(\mathbold{x}^{H_x}_i,t^{H_x}_i) = \frac{\omega_b Z_0}{N^{H_x}}  H_x^*(\mathbold{x}^{H_x}_i,t^{H_x}_i), \qquad i = 1,\dots, N^{H_x}, \label{eq:constraint_bm_1} \\
	& \frac{\omega_b Z_0}{N^{H_y}}  P_d^{H_y}(\mathbold{x}^{H_y}_i,t^{H_y}_i) = \frac{ \omega_b Z_0}{N^{H_y}}  H_y^*(\mathbold{x}^{H_y}_i,t^{H_y}_i), \qquad i = 1,\dots, N^{H_y}, \label{eq:constraint_bm_2}\\
	& \frac{\omega_b Z_0}{N^{E_z}} P_d^{E_z}(\mathbold{x}^{E_z}_i,t^{E_z}_i) = \frac{\omega_b Z_0}{N^{E_z}} E_z^*(\mathbold{x}^{E_z}_i,t^{E_z}_i), \qquad i = 1,\dots, N^{E_z}. \label{eq:constraint_bm_3}
\end{align}
Here, 
	we normalize each constraint by the total number of data available for a given variable, and $\omega_b$ is a weight that needs to be determined.
Note that pointwise values of $H_x^*$,
	$H_y^*$ and $E_z^*$ are directly given by the Hermite-Taylor method.
These constraints lead to $N^{H_x}+N^{H_y}+N^{E_z}$ additional equations that the correction functions need to satisfy.

\subsection{Boundary conditions}

In this subsection, 
    we derive constraints to enforce various boundary conditions on the surface $\Gamma$.
Consider $\Gamma$ to be parametrized by $u$ 
    and the part of the surface within a local patch $\Omega_c\cap\Gamma \times I_c = [u_\ell,u_r]\times[t_{n+1/2},t_{n+1}]$.
We define the nodes $(\mathbold{x}(u_i), t_j)$ with 
$u_i = u_\ell + i(u_r-u_\ell)/N_u$ and $t_j = t_{n+1/2} + j(t_{n+1}-t_{n+1/2})/N_u$
for $i,j = 0,\dots,N_u$ on the surface. 
In this work,
    we choose $N_u = \beta$,
    so the number of nodes $N_u+1$ discretizing $\Omega_c\cap\Gamma$ corresponds to the number of nodes on a side of the spatial domain $\Omega_{c}$ of a local patch.
Boundary conditions are then enforced on these $N_{S} = (N_u+1)^2$ nodes. 
In the following,
    we denote $(\mathbold{x}_i,t_i)$ for $i=1,\dots,N_{S}$ to be the nodes discretizing $\Omega_c\cap\Gamma \times I_c$.

Let us first consider the scaled PEC boundary condition \eqref{eq:PEC_bnd_cdn} in two space dimensions 
\begin{equation} \label{eq:2D_PEC_bnd_cdn}
    E_z = 0, \qquad c_0(n_x\mu H_x + n_y \mu H_y) = 0. 
\end{equation}
Substituting the electric field with its correction function polynomial approximation in the first equation of \eqref{eq:2D_PEC_bnd_cdn} gives
\begin{equation} \label{eq:PEC_constraint}
    E_z(\mathbold{x}_i,t_i)  \approx P_d^{E_z}(\mathbold{x}_i,t_i) = 0,  \qquad i = 1, \dots, N_S.
\end{equation}
As shown in \cite{LawAppeloHagstrom2025}, 
    we can improve the stability of the overall method, 
    that is, 
    when the CFM is applied to the Hermite-Taylor method, 
    by constraining certain spatial derivatives of the correction functions on the surface. 
In this situation, 
    we consider time derivatives of boundary conditions, 
    which are converted into spatial derivatives by using the governing equations.
Although increasing the maximum order 
    of the time derivatives to be constrained improves the stability of the overall method, 
    it usually also increases the condition number of the resulting matrix involved in the minimization procedure.
Considering smooth solutions of Maxwell's equations \eqref{eq:Maxwell_equations}, 
    we have 
\begin{align}
    & \partial_t^{j} \mathbold{E} = \left.\Bigg\{ 
    \begin{aligned}
        \nabla^{j-1}(\nabla\times \mathbold{H})/(\epsilon^{\theta}\,\mu^{\theta-1}), \qquad &\mbox{if } j \mbox{ odd}, \\
        -\nabla^{j-2}(\nabla\times\nabla\times \mathbold{E})/(\epsilon\,\mu)^{j/2}, \qquad &\mbox{otherwise},
    \end{aligned}\right. \label{eq:convert_dtE} \\
    & \partial_t^{j} \mathbold{H} = \left.\Bigg\{ 
    \begin{aligned}
        -\nabla^{j-1}(\nabla\times \mathbold{E})/(\epsilon^{\theta-1}\,\mu^{\theta}), \qquad &\mbox{if } j \mbox{ odd}, \\
        -\nabla^{j-2}(\nabla\times\nabla\times \mathbold{H})/(\epsilon\,\mu)^{j/2}, \qquad &\mbox{otherwise},
    \end{aligned}\right. \label{eq:convert_dtH}
\end{align}
with $\theta = (j+1)/2$.
In addition to constraints \eqref{eq:PEC_constraint}, 
    we also consider   
\begin{equation}
\bigg(\frac{L_0}{c_0}\bigg)^j\partial_t^j P_d^{E_z} (\mathbold{x}_i,t_i) = 0, \qquad i = 1,\dots, N_S, \qquad  j=1,\dots,N_d,
\end{equation}
    where all the time derivatives $\partial_t^j P_d^{E_z}$ for  
    $j=1,\dots,N_d$ 
    are converted into spatial derivatives using \eqref{eq:convert_dtE}.
Even though we do not explicitly write this conversion in every case for legibility,  
    we always consider that all the time derivatives used to enforce surface conditions are converted into spatial derivatives throughout the remainder of this work.
Repeating the procedure for the second equation of \eqref{eq:2D_PEC_bnd_cdn}, 
    the constraints for PEC boundary condition \eqref{eq:2D_PEC_bnd_cdn} are then 
\begin{align} 
& \frac{\omega_s}{N_S}\bigg(\frac{L_0}{c_0}\bigg)^j\partial_t^j P_d^{E_z} (\mathbold{x}_i,t_i) = 0, \label{eq:constraints_PEC_1} \\
& \frac{\omega_s}{N_S}c_0\bigg(\frac{L_0}{c_0}\bigg)^j(n_x  \mu \partial_t^j P_d^{H_x}(\mathbold{x}_i,t_i) + n_y  \mu \partial_t^j P_d^{H_y}(\mathbold{x}_i,t_i) )=0, \label{eq:constraints_PEC_2}
\end{align}
for $i = 1,\dots, N_S$ and $j=0,\dots,N_d$. 
Here $\omega_s$ is a given weight and $N_d \leq d$.
We also normalize each equation by the total number of points $N_s$.

In two space dimensions, 
    PMC boundary condition \eqref{eq:PMC_bnd_cdn} with the appropriate scaling is simplified to  
\begin{equation} \label{eq:2D_PMC_bnd_cdn}
 Z_0(n_x H_y - n_y H_x) = 0, 
\end{equation} 
    while impedance boundary condition \eqref{eq:impedance_bnd_cdn} becomes
\begin{equation} \label{eq:2D_impedance_bnd_cdn}
    E_z - Z (n_y H_x - n_x H_y) = 0.
\end{equation}
Repeating the same procedure for PMC and impedance boundary conditions leads respectively to 
\begin{equation}
    \frac{\omega_s}{N_S}Z_0\bigg(\frac{L_0}{c_0}\bigg)^j \big(n_x\partial_t^j P_d^{H_y}(\mathbold{x}_i,t_i) - n_y \partial_t^j P_d^{H_x}(\mathbold{x}_i,t_i)\big) =0, \qquad i = 1,\dots, N_S, \qquad j=0,\dots,N_d, 
\end{equation}
and 
\begin{equation}
    \frac{\omega_s}{N_S}\bigg(\frac{L_0}{c_0}\bigg)^j \big(\partial_t^j P_d^{E_z}(\mathbold{x}_i,t_i) - Z (n_y\partial_t^j P_d^{H_x}(\mathbold{x}_i,t_i) + n_x \partial_t^j P_d^{H_y}(\mathbold{x}_i,t_i))\big) =0, \qquad i = 1,\dots, N_S, \qquad j=0,\dots,N_d.
\end{equation}
In all cases, 
    we have$(N_d+1)N_S$ additional  constraints to enforce. 

\subsection{Solving the linear system for the correction function polynomial coefficients} \label{sec:solving_linear_system_bnd_cdns}

The above constraints lead to a system of equations of the form $M \mathbold{c} = \mathbold{b}$ with $M = [\mathcal{G}; \mathcal{B}; \mathcal{S}]$ and $\mathbold{b} = [0; \mathbold{b}_{\mathcal{B}} ; 0]$.
Here, 
    the matrices $\mathcal{G}$,
    $\mathcal{B}$ and $\mathcal{S}$ are coming from enforcing the governing equations, 
    matching the base method's numerical solution and enforcing the boundary conditions, 
    respectively. 
Noticing that we have a total of $N_{eq} = 3d^3+7d^2+7d-2 + N^{H_x} + N^{H_y} + N^{E_z} + (N_d+1)N_s$ equations for $N_c = 2d^3+6d^2+7d+3$ polynomial coefficients to determine, 
    the linear system of equations is then overdetermined.
Actually, 
    for $d>1$, 
    the number of rows of $\mathcal{G}$ is already greater than the number of components of $\mathbold{c}$.

The matrix $\mathcal{G}$ is rank deficient since constraints \eqref{eq:constraints_Faraday_1}-\eqref{eq:constraints_div_last}, 
    enforcing the governing equations, 
    do not involve  coefficients $c_{0,0,0}^{H_x}$, 
$c_{0,0,0}^{H_y}$ and $c_{0,0,0}^{E_z}$.
For all boundary conditions and their time derivatives, 
    we do not necessarily enforce all components of the electromagnetic fields on the surface $\Gamma$, 
    so the matrix $\mathcal{S}$ is also rank deficient. 
However, 
    noticing that the matrix $\mathcal{B}$ is actually a rectangular generalized Vandermonde matrix, 
    we obtain a full column rank matrix $M$ for a sufficiently large number of well-distributed BM nodes,
    taken among the Cartesian primal and dual meshes used by the Hermite-Taylor method.
In this situation, 
    the least squares problem 
\begin{equation} \label{eq:min_problem}
    \displaystyle\min_{\mathbold{c} \in \mathbb{R}^{N_c}}\|M\mathbold{c} - \mathbold{b}\|_2
\end{equation}
    has a unique solution.
We therefore solve the linear system of equations in the least squares sense 
    using a pivoted QR factorization. 
Since the surface is independent of time, 
    the construction of local patches and the QR factorization of the associated matrices $M$ are computed in a pre-computation step. 
For each half time step and for each patch, 
    we solve problem \eqref{eq:min_problem} and evaluate the correction function polynomial approximations and their derivatives through order $m$ at the CF nodes associated with the local patch.

\subsection{Coupling the Hermite-Taylor method and the DCFM}

In this subsection, 
    we describe the algorithm to couple the Hermite-Taylor method with the DCFM over one time step.
Assume that approximations of all the variables and their derivatives through order $m$ are available on the primal mesh at time $t_n$.
The algorithm of the Hermite-Taylor discrete correction function method to evolve these data on the primal mesh at time $t_{n+1}$ is 
\begin{itemize}
    \item[1.] Solve for the dual BM nodes at time $t_{n+1/2}$ using the Hermite-Taylor method;
    \item[2.] Solve for the dual CF nodes at time $t_{n+1/2}$ using the DCFM. 
    The time domain of local patches is $I_c = [t_n, t_{n+1/2}]$, 
        and the DCFM involves certain primal BM nodes at $t_n$ and dual BM nodes at $t_{n+1/2}$;
    \item[3.] Solve for the primal BM nodes at time $t_{n+1}$ using the Hermite-Taylor method;
    \item[4. ]Solve for the primal CF nodes at time $t_{n+1}$ using the DCFM. 
    The time domain of local patches is $I_c = [t_{n+1/2}, t_{n+1}]$, 
        and the DCFM involves certain dual BM nodes at $t_{n+1/2}$ and primal BM nodes at $t_{n+1}$.
\end{itemize}
Note that,
    at each half time step, 
    minimization problem \eqref{eq:min_problem} for each local patch can be solved independently, 
    making it suitable for parallel implementations \cite{Abraham2017}.

\subsection{Extension to variable coefficient problems} \label{sec:var_coeff_problems}

Let us now consider that the magnetic permeability and the electric permittivity depend on the spatial coordinates and are such that $\mu(\mathbold{x}),\epsilon(\mathbold{x}) > 0$.
In this situation, 
    we consider the constitutive laws 
\begin{equation} \label{eq:constitutive_laws_var_coeff}
    \mathbold{B} = \mu \mathbold{H}, \qquad \mathbold{D} = \epsilon \mathbold{E},
\end{equation}
    where $\mathbold{B}$ and $\mathbold{D}$ are respectively the magnetic induction and the electric displacement,
    and we rewrite Maxwell's equations as
\begin{equation} \label{eq:Maxwell_equations_B_D}
		\partial_t \mathbold{B} + \nabla\times (\mathbold{D}/\epsilon) = 0, 
        \qquad
		\partial_t \mathbold{D} - \nabla\times (\mathbold{B}/\mu) = 0, \qquad
		\nabla\cdot\mathbold{D} = 0, \qquad
		\nabla\cdot\mathbold{B} = 0.
\end{equation}
Note that the electromagnetic fields we now seek are divergence-free,
    and therefore the discrete correction function method introduced previously can be used with few modifications, 
    as explained below.

Let us focus on the governing equations constraints.
Assuming sufficiently smooth coefficients $1/\mu$ and $1/\epsilon$, 
    we have 
\begin{equation}
    \frac{1}{\mu(\mathbold{x})} \approx P_d^{1/\mu}(\mathbold{x}) = \sum_{k=0}^d\sum_{\ell=0}^d c_{k,\ell}^{1/\mu}\xi^k\eta^\ell, \qquad \frac{1}{\epsilon(\mathbold{x})} \approx P_d^{1/\epsilon}(\mathbold{x}) = \sum_{k=0}^d\sum_{\ell=0}^d c_{k,\ell}^{1/\epsilon}\xi^k\eta^\ell.
\end{equation}
Note that coefficients $c_{k,\ell}^{1/\mu}$ and $c_{k,\ell}^{1/\epsilon}$ correspond to scaled spatial derivatives of $1/\mu$ and $1/\epsilon$, 
    and therefore are known. 
Substituting the electromagnetic fields and the variable coefficients by their polynomial approximations in the first equation of Faraday's law in two space dimensions,
    we obtain 
\begin{equation} \label{eq:Faraday_law_1_var_coeff}
    L_0(\partial_t P_d^{B_x} + \partial_y (P_d^{1/\epsilon}P_d^{D_z})) = 0.
\end{equation}
Noticing $E_z \approx P_d^{1/\epsilon} P_d^{D_z}$ and 
truncating this polynomial product to obtain a polynomial of degree at most $d$ in each independent variable leads to 
\begin{equation}
E_z \approx \sum_{k=0}^d \sum_{\ell=0}^d\sum_{s=0}^{d} c_{k,\ell,s}^{E_z} \xi^k\eta^\ell\zeta^s,
\end{equation}
with 
\begin{equation} \label{eq:coeff_Ez_var_coeff}
c_{k,\ell,s}^{E_z} = \sum_{i=0}^k\sum_{j=0}^\ell c_{i,j}^{1/\epsilon} c_{k-i,\ell-j,s}^{D_z}.
\end{equation}
Matching the polynomial coefficients in \eqref{eq:Faraday_law_1_var_coeff} provides the following equations for the coefficients of $B_x$ and $D_z$
\begin{align}
   & \frac{L_0}{\Delta t}(s+1) c_{k,\ell,s+1}^{B_x}  + \frac{L_0}{\Delta y}(\ell+1)c_{k,\ell+1,s}^{E_z} = 0, \qquad k=0,\dots,d, \qquad \ell,s = 0,\dots,d-1, \label{eq:constraints_Faraday_1_var_coeff} \\
   & \frac{L_0}{\Delta t}   (s+1) c_{k,d,s+1}^{B_x} = 0, \qquad k = 0,\dots,d, \qquad s = 0,\dots, d-1, \label{eq:constraints_Faraday_2_var_coeff}\\
   & \frac{L_0}{\Delta y} (\ell+1) c_{k,\ell+1,d}^{E_z} = 0, \qquad k = 0,\dots,d, \qquad \ell = 0,\dots, d-1, \label{eq:constraints_Faraday_3_var_coeff}
\end{align}
    with coefficients $c_{k,\ell,s}^{E_z}$ given by \eqref{eq:coeff_Ez_var_coeff}. 
The remaining equations for the correction function polynomial coefficients coming from the second equation of Faraday's law and Amp\`ere-Maxwell's law are obtained using a similar procedure.
In this situation, 
    we cannot leverage the divergence-free constraint $\nabla\cdot\mathbold{B}=0$ to simplify Amp\`ere-Maxwell's law as it is done for the constant coefficient case, 
    that is, 
    when we simplify equations
    \eqref{eq:constraints_Ampere_Maxwell_6} and 
    \eqref{eq:constraints_Ampere_Maxwell_7} to equations
    \eqref{eq:constraints_Ampere_Maxwell_6_modify} and 
    \eqref{eq:constraints_Ampere_Maxwell_7_modify} in subsection \ref{sec:gov_eq_constraints}.

Constraints involving the base method's solution require computing pointwise values of $\mathbold{B}^*$ and $\mathbold{D}^*$ from $\mathbold{H}^*$ and $\mathbold{E}^*$, 
    obtained using the Hermite-Taylor method.
This can be easily done using constitutive laws \eqref{eq:constitutive_laws_var_coeff}. 
Regarding boundary conditions involving  spatially variable coefficients, 
    such as PEC boundary condition, 
    we directly evaluate the physical parameters at nodes $\mathbold{x}_i$ discretizing the surface $\Gamma$.
For additional constraints involving high-order time derivatives 
    of a boundary condition,
    knowing  the spatial derivatives of $1/\mu(\mathbold{x})$ and $1/\epsilon(\mathbold{x})$,
    we convert time derivatives of the electromagnetic fields into spatial derivatives by leveraging Leibniz's rule.

\subsection{Extension to interface and generalized sheet transition conditions}

The generalization to interface and metasurface problems requires considering two sets of correction functions, 
    one for each subdomain.
We then define the polynomial correction function $P_d^{f,\pm}$ approximating the variable $f$ for the subdomain $\Omega^{\pm}$.
In this situation, 
    we also have two sets of constraints $\mathcal{G}^\pm \mathbold{c}^\pm = 0$ for enforcing the governing equations in each subdomain. 
Correction functions are also required to match the base method's numerical solution at nodes within their respective subdomains in a local patch.
We then have $\mathcal{B}^\pm \mathbold{c}^\pm = \mathbold{b}_{\mathcal{B}}^\pm$.

Let us first consider conventional interface conditions \eqref{eq:interface_cdns} with the appropriate scaling in two space dimensions
\begin{equation} \label{eq:interface_cdns_2d}
     Z_0(\tilde{n}_x \llbracket H_y \rrbracket - \tilde{n}_y \llbracket H_x \rrbracket) = 0, \qquad c_0(\tilde{n}_x \llbracket \mu H_x \rrbracket + \tilde{n}_y \llbracket \mu H_y \rrbracket) = 0, \qquad \llbracket E_z \rrbracket = 0.
\end{equation}
In this situation, 
    the procedure to derive the constraints to enforce interface conditions is the same as for boundary conditions since there are neither time nor spatial derivatives involved in these conditions.
Applying time derivatives to interface  conditions \eqref{eq:interface_cdns_2d}, 
    substituting components of the electromagnetic fields with their polynomial approximations and enforcing the resulting constraints at the nodes discretizing the part of the surface $\tilde{\Gamma}\times I$ included in a local patch leads to 
\begin{align}
& \frac{\omega_s}{N_S}Z_0\bigg(\frac{L_0}{c_0}\bigg)^j(\tilde{n}_x \llbracket \partial_t^j P_d^{H_y}(\mathbold{x}_i,t_i) \rrbracket - \tilde{n}_y \llbracket \partial_t^j P_d^{H_x}(\mathbold{x}_i,t_i) \rrbracket) = 0, \\
& \frac{\omega_s}{N_S}c_0\bigg(\frac{L_0}{c_0}\bigg)^j(\tilde{n}_x \llbracket \mu \partial_t^j P_d^{H_x}(\mathbold{x}_i,t_i) \rrbracket + \tilde{n}_y \llbracket \mu \partial_t^j P_d^{H_y}(\mathbold{x}_i,t_i) \rrbracket)=0, \\ 
& \frac{\omega_s}{N_S}\bigg(\frac{L_0}{c_0}\bigg)^j\llbracket \partial_t^j P_d^{E_z} (\mathbold{x}_i,t_i) \rrbracket = 0, 
\end{align}
for $i=1,\dots,N_s$ and $j=0,\dots,N_d$.

Let us now consider GSTC model \eqref{eq:GSTC_time_domain},
    involving time and spatial derivatives of the electromagnetic fields, 
    in two space dimensions.
Assuming a metasurface to be flat with the normal  
    $\tilde{\mathbold{n}} = [1;0]$ leads to
\begin{equation} \label{eq:GSTC_flat}
Z_0(\llbracket H_y \rrbracket - \epsilon_0 \chi_{ee}^{zz} \partial_t \clbracket E_z \crbracket) = 0, \qquad 
c_0(\llbracket \mu H_x \rrbracket + \mu_0  \chi_{mm}^{yy} \partial_y \clbracket H_y \crbracket) =0, \qquad 
\llbracket E_z \rrbracket - \mu_0 \chi_{mm}^{yy} \partial_t \clbracket H_y \crbracket = 0.
\end{equation}
In this situation,
    we convert time derivatives 
    $\partial_t \clbracket H_y \crbracket$ and $\partial_t \clbracket E_z \crbracket$ into spatial derivatives using equations \eqref{eq:convert_dtE} and \eqref{eq:convert_dtH}.
We then follow the same procedure as before to obtain the following constraints 
\begin{align}
\frac{\omega_s}{N_S}Z_0\bigg(\frac{L_0}{c_0}\bigg)^j(\llbracket \partial_t^j P_d^{H_y}(\mathbold{x}_i,t_i) \rrbracket - \epsilon_0 \chi_{ee}^{zz}  \clbracket \partial_t^{j+1}  P_d^{E_z}(\mathbold{x}_i,t_i) \crbracket) = 0, \label{eq:constraints_GSTC_1} \\
\frac{\omega_s}{N_S}c_0\bigg(\frac{L_0}{c_0}\bigg)^j(\llbracket \mu \partial_t^j P_d^{H_x}(\mathbold{x}_i,t_i) \rrbracket + \mu_0  \chi_{mm}^{yy} \clbracket \partial_t^j \partial_y P_d^{H_y}(\mathbold{x}_i,t_i) \crbracket) =0, \label{eq:constraints_GSTC_2} \\
\frac{\omega_s}{N_S}\bigg(\frac{L_0}{c_0}\bigg)^j (\llbracket \partial_t^j P_d^{E_z}(\mathbold{x}_i,t_i) \rrbracket - \mu_0 \chi_{mm}^{yy} \clbracket \partial_t^{j+1} P_d^{H_y}(\mathbold{x}_i,t_i) \crbracket) = 0, \label{eq:constraints_GSTC_3}
\end{align}
    for $i=1,\dots,N_s$ and $j=0,\dots,N_d$.
Since GSTC model \eqref{eq:GSTC_time_domain} involves first-order derivatives,
    we have $N_d < d$ for all the terms in the constraints to be considered. 

For interface and generalized sheet transition conditions, 
    we then have the following system $M\mathbold{c} = \mathbold{b}$ with
\begin{equation}
    M = 
    \begin{bmatrix}
    \mathcal{G} \\
    \mathcal{B} \\
    \mathcal{S}
    \end{bmatrix},
    \qquad 
    \mathcal{G} = 
    \begin{bmatrix}
        \mathcal{G}^+ & 0 \\
        0 & \mathcal{G}^- 
    \end{bmatrix},
    \qquad 
    \mathcal{B} = 
    \begin{bmatrix}
        \mathcal{B}^+ & 0 \\
        0 & \mathcal{B}^- 
    \end{bmatrix},
    \qquad
    \mathbold{c} = 
    \begin{bmatrix}
    \mathbold{c}^+ \\ 
    \mathbold{c}^-
    \end{bmatrix},
    \qquad 
    \mathbold{b} =
    \begin{bmatrix}
    0 \\
    0 \\
    \mathbold{b}_{\mathcal{B}}^+ \\ 
    \mathbold{b}_{\mathcal{B}}^- \\
    0 
    \end{bmatrix}.
\end{equation}
This system is also overdetermined and the least squares problem is solved using a pivoted QR factorization, 
    as described in subsection~\ref{sec:solving_linear_system_bnd_cdns}.

\subsection{Determination of the weights}

The choice of the weights $\omega_b$ and $\omega_s$ is determined using a local truncation error analysis. 
Assume the base method provides a numerical solution of order $k$ and that polynomials approximating the correction functions have an accuracy of $\mathcal{O}(h^{d+1})$.
We first consider the constraints requiring the correction functions to match the base method's numerical solution.
Considering constraint \eqref{eq:constraint_bm_1}, 
    we then have 
\begin{equation}
\frac{\omega_b Z_0}{N^{H_x}} (H_x -P_d^{H_x}) -\frac{\omega_b Z_0}{N^{H_x}}  (H_x -H_x^*) = \mathcal{O}(h^k) + \mathcal{O}(h^{d+1}).
\end{equation}
So the local truncation error is of order $\min\{k,d+1\}$. 
Constraints \eqref{eq:constraint_bm_2} and \eqref{eq:constraint_bm_3} lead to a local truncation error of the same order.
Note that we choose $0 < \omega_b < 1$ to emphasize the condition to be enforced on the surface.

Let us consider the constraints for PEC boundary condition \eqref{eq:2D_PEC_bnd_cdn}.
We obtain 
\begin{equation}
\frac{\omega_s}{N_S}\bigg(\frac{L_0}{c_0}\bigg)^j\partial_t^j (E_z-P_d^{E_z})= \mathcal{O}(\omega_sL_0^j h^{d+1-j}),  \qquad   j=0,\dots,N_d.
\end{equation}
Since $L_0 = \mathcal{O}(h)$, 
    we have a local truncation error of order $d+1$ for all $j$.
Repeating the analysis for PMC boundary condition, 
    impedance boundary condition, 
    and interface conditions 
    also results in a local truncation error of order $d+1$.
By choosing $d \geq k - 1$ to match the accuracy of the base method, 
    all the terms scale as $\mathcal{O}(h^k)$ as the mesh size diminishes. 
We therefore set $\omega_s=1$ for these conditions.

Consider now constraints \eqref{eq:constraints_GSTC_1} coming from GSTC model \eqref{eq:GSTC_flat}.
We have 
\begin{equation} \label{eq:GSTC_eq_1_truncation_err}
\frac{\omega_s}{N_S}Z_0\bigg(\frac{L_0}{c_0}\bigg)^j(\llbracket \partial_t^j (H_y-P_d^{H_y}) \rrbracket - \epsilon_0 \chi_{ee}^{zz}  \clbracket \partial_t^{j+1} (E_z - P_d^{E_z}) \crbracket) = \mathcal{O}(\omega_sL_0^j(h^{d+1-j}+h^{d-j})).
\end{equation}
Providing $L_0 = \mathcal{O}(h)$, 
    the local error \eqref{eq:GSTC_eq_1_truncation_err} is of order $d$.
Constraints \eqref{eq:constraints_GSTC_2} and \eqref{eq:constraints_GSTC_3} also have a local truncation error of order $d$ because of first-order derivatives involved in GSTC model \eqref{eq:GSTC_flat}. 
We therefore set $\omega_s=L_0$ to have all the terms scaled as $\mathcal{O}(h^k)$ when $d\geq k-1$.

\section{Numerical examples} \label{sec:num_examples}

In this section, 
    we investigate the properties of the proposed Hermite-Taylor discrete correction function method by performing long-time simulations and convergence studies in two space dimensions. 
In all the numerical examples, 
    the computational domain $\Omega$ is discretized using a Cartesian mesh with the mesh size $h = \Delta x = \Delta y$.
We also set $\omega_b = 1/2$. 
We consider the surface to be either a circle or a smooth three-star shape, 
    described by the level set function
\begin{equation} \label{eq:level_set_fct}
\phi(x,y) = (x-x_0)^2 + (y-y_0)^2 - r^2(x,y), \qquad r(x,y) = r_0 +  a\sin(b\theta(x,y)),
\end{equation}
with $a=1/10$, 
    $b=3$, 
    $r_0 = 0.3$ and $\theta(x,y)$ is the angle between the $x$-axis and the vector $[x-x_0; y-y_0]$.

The error of the electromagnetic fields at the final time is computed using the usual 2-norm. 
Regarding the divergence-free constraint, 
    we use the $L_2$-norm of the divergence of the Hermite-Taylor and the correction function polynomials at the final time, 
    as in \cite{LawAppeloHagstrom2025}. 
Considering a domain divided into two subdomains $\Omega^{\pm}$. 
The $L_2$-norm of the divergence of the magnetic field at the final time is computed using 
\begin{equation} \label{eq:computation_divergence_H}
\|\nabla\cdot(\mu\mathbold{H}_h)\|_2 = \|\nabla\cdot(\mu\mathbold{H}^*)\|_2 + \sum_{j=1}^{N_{c}} \Big(\|\nabla\cdot(\mu\mathbold{H}^+_{h,j})\|_2 + 
    \|\nabla\cdot(\mu\mathbold{H}^-_{h,j})\|_2\Big).
\end{equation}
Here $N_{c}$ is the number of local patches, 
    $\mathbold{H}^*$ is the Hermite-Taylor approximation of the magnetic field, 
    and 
    $\mathbold{H}_{h,j}^\pm$ is the approximation of the magnetic field associated with the $j$th local patch in subdomain $\Omega^{\pm}$.

\subsection{Boundary conditions}

We consider the physical domain $\Omega_p$, 
    for which its boundary $\Gamma$ is described by the level set function \eqref{eq:level_set_fct} with $x_0=y_0=0.5$.
The domain $\Omega_p$ is immersed in the computational domain $\Omega = [0,1]\times[0,1]$. 
The physical parameters are $\mu = 1$ and $\epsilon =1$. 
The initial condition and boundary conditions on $\Gamma$ are chosen so that the standing mode solution 
\begin{equation} \label{eq:sol_bnd_cdn_cst_coeff}
    \begin{aligned}
      &   H_x = -\sin(\omega \pi x) \cos(\omega \pi y) \sin(\sqrt{2}\omega\pi t)/\sqrt{2}, \\ 
      &   H_y = \cos(\omega \pi x) \sin(\omega \pi y) \sin(\sqrt{2}\omega\pi t)/\sqrt{2}, \\ 
      &  E_z = \sin(\omega\pi x) \sin(\omega \pi y) \cos(\sqrt{2}\omega\pi t), 
    \end{aligned}
\end{equation}
    with $\omega = 5$ is satisfied in the physical domain $\Omega_p$. 
Note that source terms in boundary conditions \eqref{eq:2D_PEC_bnd_cdn}, 
    \eqref{eq:2D_PMC_bnd_cdn} and \eqref{eq:2D_impedance_bnd_cdn} are necessary.

\begin{figure}   
	\centering
	\begin{adjustbox}{max width=0.9\textwidth,center}
		\includegraphics[width=\linewidth,trim={0cm 0cm 1.75cm 0cm},clip]{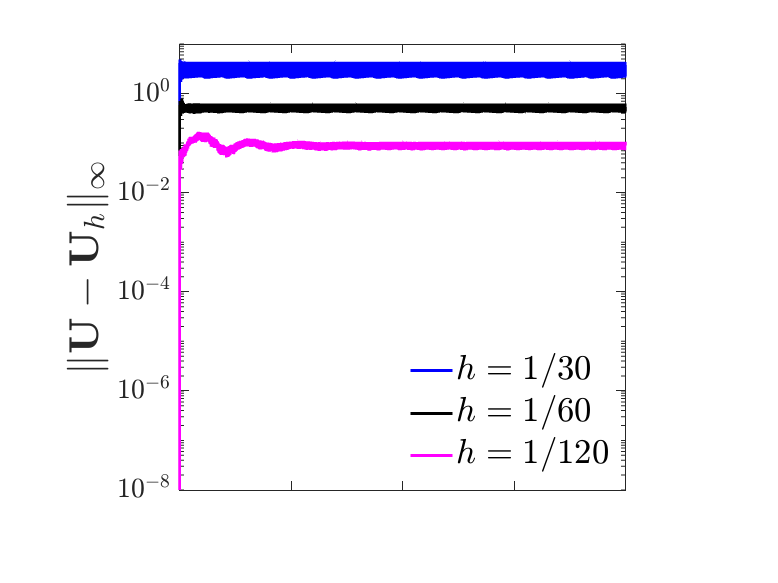} \hspace{-18.0pt}
		\includegraphics[width=\linewidth,trim={0cm 0cm 1.75cm 0cm},clip]{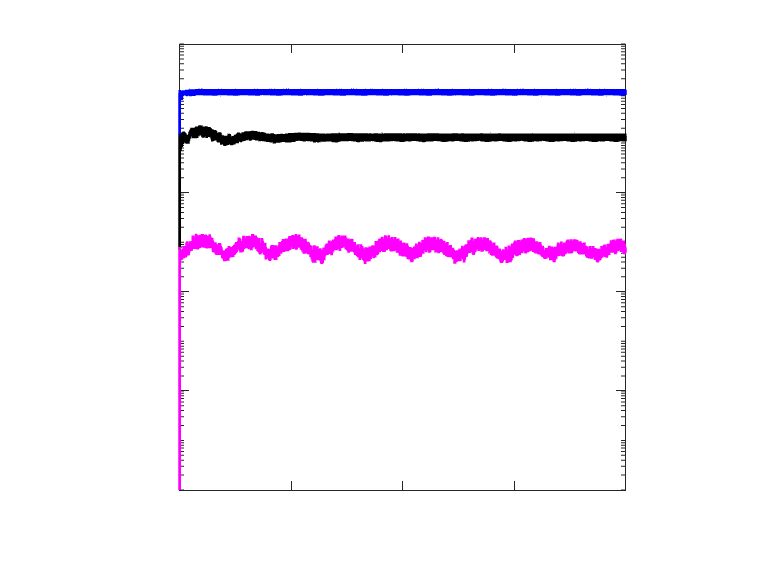}\hspace{-18pt}
		\includegraphics[width=\linewidth,trim={0cm 0cm 1.75cm 0cm},clip]{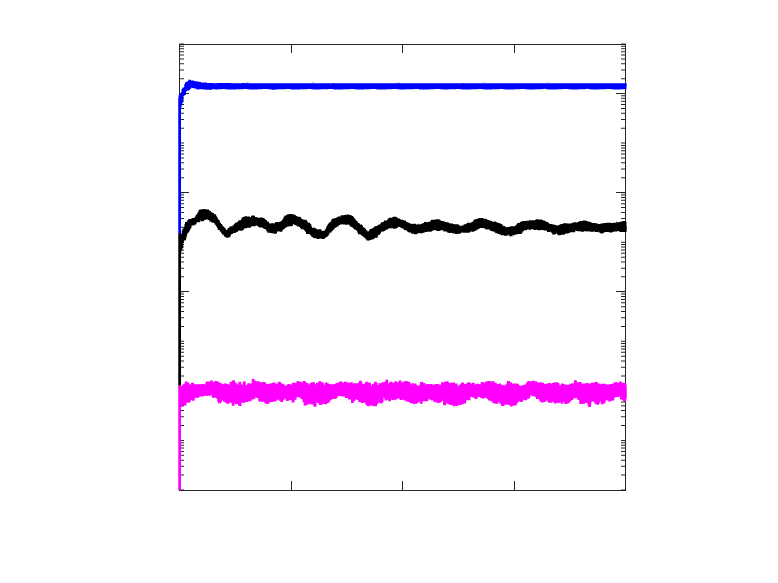}\hspace{-18pt}
	\end{adjustbox} 
	\begin{adjustbox}{max width=0.9\textwidth,center}
		\includegraphics[width=\linewidth,trim={0cm 0cm 1.75cm 0cm},clip]{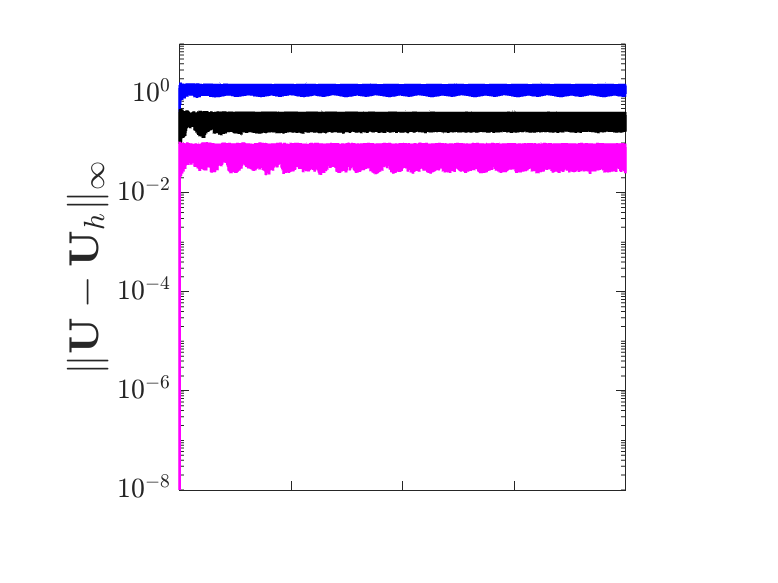} \hspace{-18.0pt}
		\includegraphics[width=\linewidth,trim={0cm 0cm 1.75cm 0cm},clip]{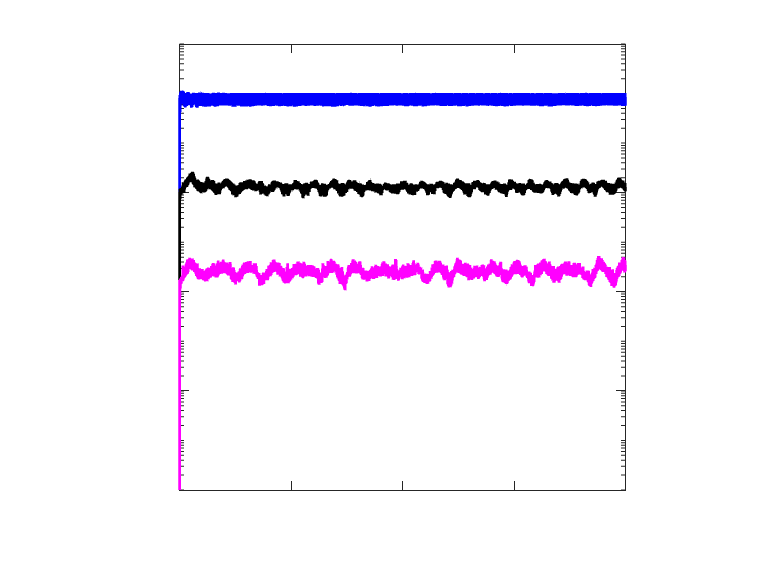}\hspace{-18pt}
		\includegraphics[width=\linewidth,trim={0cm 0cm 1.75cm 0cm},clip]{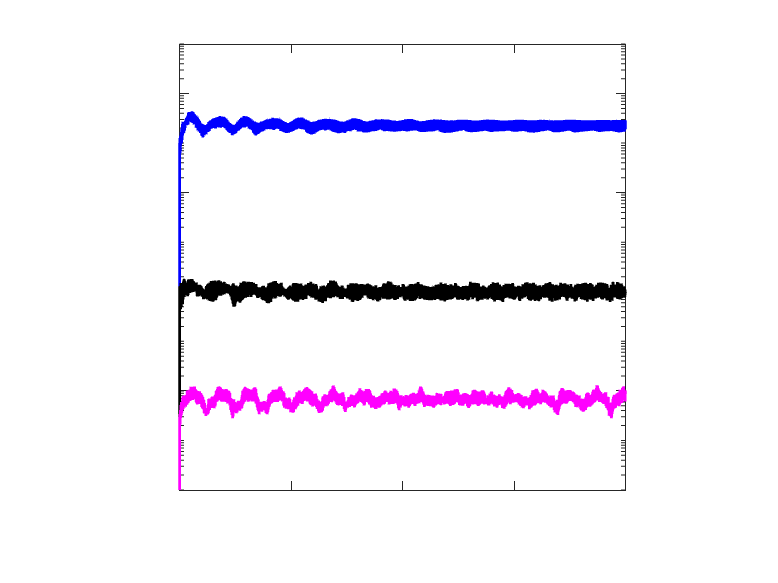}\hspace{-18pt}
	\end{adjustbox} 
	\begin{adjustbox}{max width=0.9\textwidth,center}
		\includegraphics[width=\linewidth,trim={0cm 0cm 1.75cm 0cm},clip]{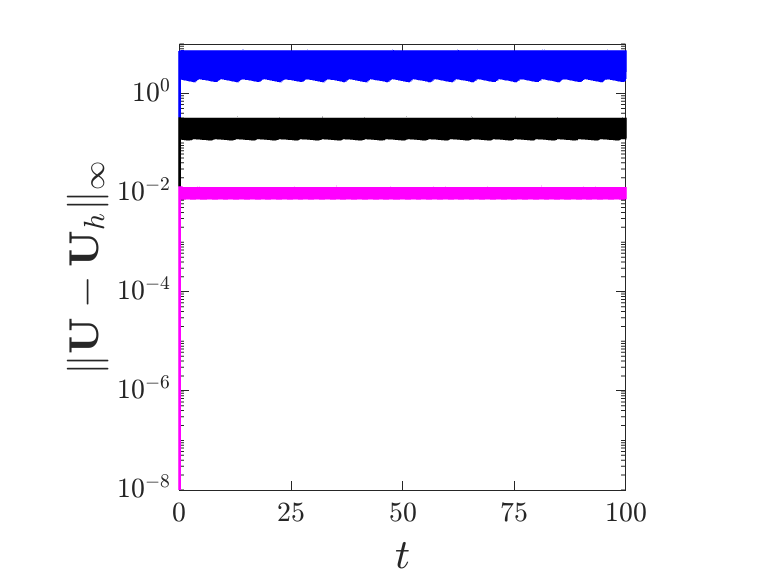} \hspace{-18.0pt}
		\includegraphics[width=\linewidth,trim={0cm 0cm 1.75cm 0cm},clip]{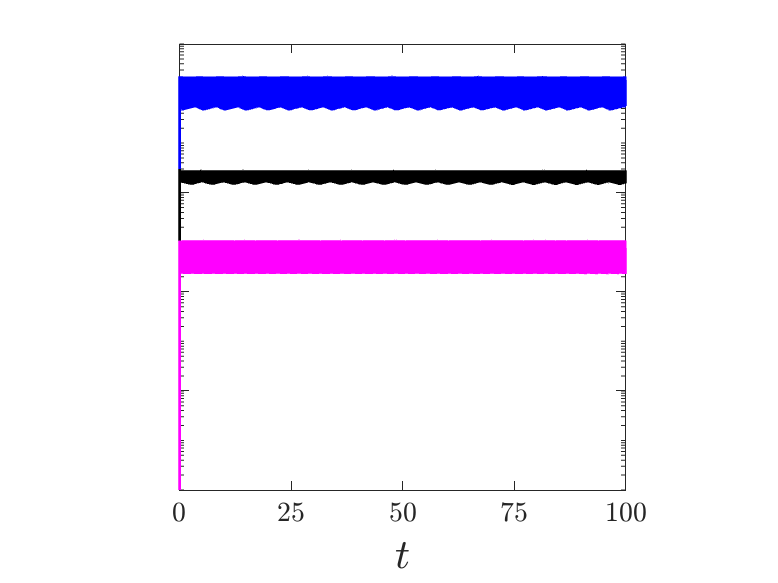}\hspace{-18pt}
		\includegraphics[width=\linewidth,trim={0cm 0cm 1.75cm 0cm},clip]{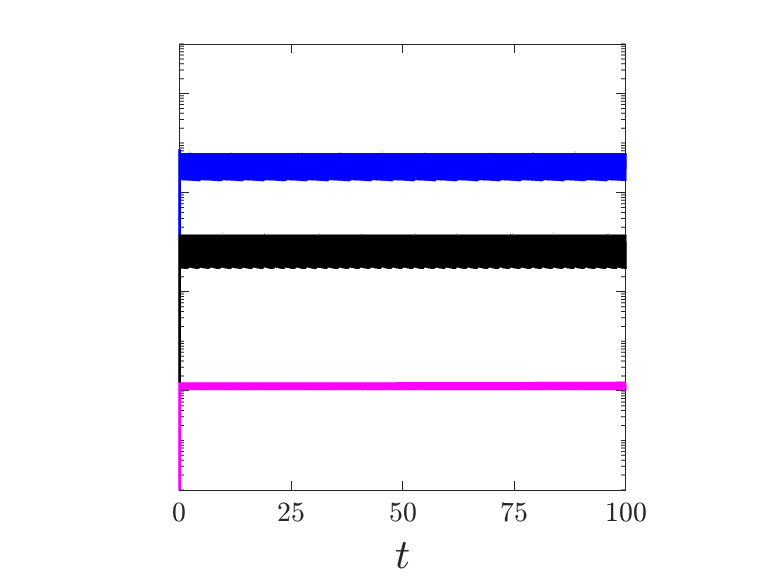}\hspace{-18pt}
	\end{adjustbox} 
       \caption{Maximum norm of the error as a function of time for different mesh sizes, boundary conditions and values of $m$. The rows correspond to problems using  different boundary conditions: PEC, PMC and impedance from top to bottom, 
       respectively. The columns are for different m: 1 to 3 from left to right. 
       Here $\mathbold{U} = [H_x;H_y;E_z]$.}
\label{fig:long_time_simulations_bnd_cdns}
\end{figure}

We first perform long-time simulations to numerically investigate the stability of the proposed method for $m=1-3$. 
The time interval is $I = [0,100]$.
In this situation, 
    a wave travels 250 wavelengths and there 
    are approximately 112.5 periods. 
We consider three different meshes with $h\in\{1/30, 1/60, 1/120\}$.
The correction functions are sought in the space of divergence-free polynomials of degree at most $d=2m$, 
    and $N_d = 2m$.
For PEC and PMC boundary conditions, 
    we choose the CFL constant to be $0.5$ for all values of $m$. 
For impedance boundary condition, 
    the CFL constant is $0.5$ for $m=1-2$ and $0.25$ for $m=3$. 
Fig.~\ref{fig:long_time_simulations_bnd_cdns} illustrates the error in maximum norm as a function of time. 
The results suggest that the method is stable.

The maximum condition number of  matrices $M$ coming from the three different meshes is around $10^4$,
    $10^6$, 
    $10^8$ for $m=1-3$, 
    respectively,
    using PEC boundary condition \eqref{eq:2D_PEC_bnd_cdn}.
As mentioned before, 
    the DCFM avoids the normal form of the least squares problem,
    providing better-conditioned matrices when compared to the CFM in \cite{LawAppeloHagstrom2025}. 
For PMC boundary condition \eqref{eq:2D_PMC_bnd_cdn}, 
    the maximum condition number is approximately $10^4$, 
    $10^7$ 
    and $10^{10}$, 
    while for impedance boundary condition \eqref{eq:2D_impedance_bnd_cdn}, 
    we have around $10^5$,
    $10^8$ 
    and $10^{11}$,
    for $m=1-3$, 
    respectively.

The DCFM becomes unstable for $m=4$, 
    corresponding to a ninth-order method,
    when PEC boundary condition is considered.
Note that we did not investigate CFL constants lower than $0.1$,
    as we aim to retain the remarkable stability properties of the Hermite-Taylor method while enforcing boundary conditions. 
A condition number of around $10^{12}$ is obtained with the CFL constant of 0.1.

\begin{figure}
\centering
\begin{subfigure}{1.0\textwidth}
    \centering
    \begin{adjustbox}{max width=0.9\textwidth,center}
    \includegraphics[width=\linewidth,trim={0cm 0cm 1.75cm 0cm},clip]{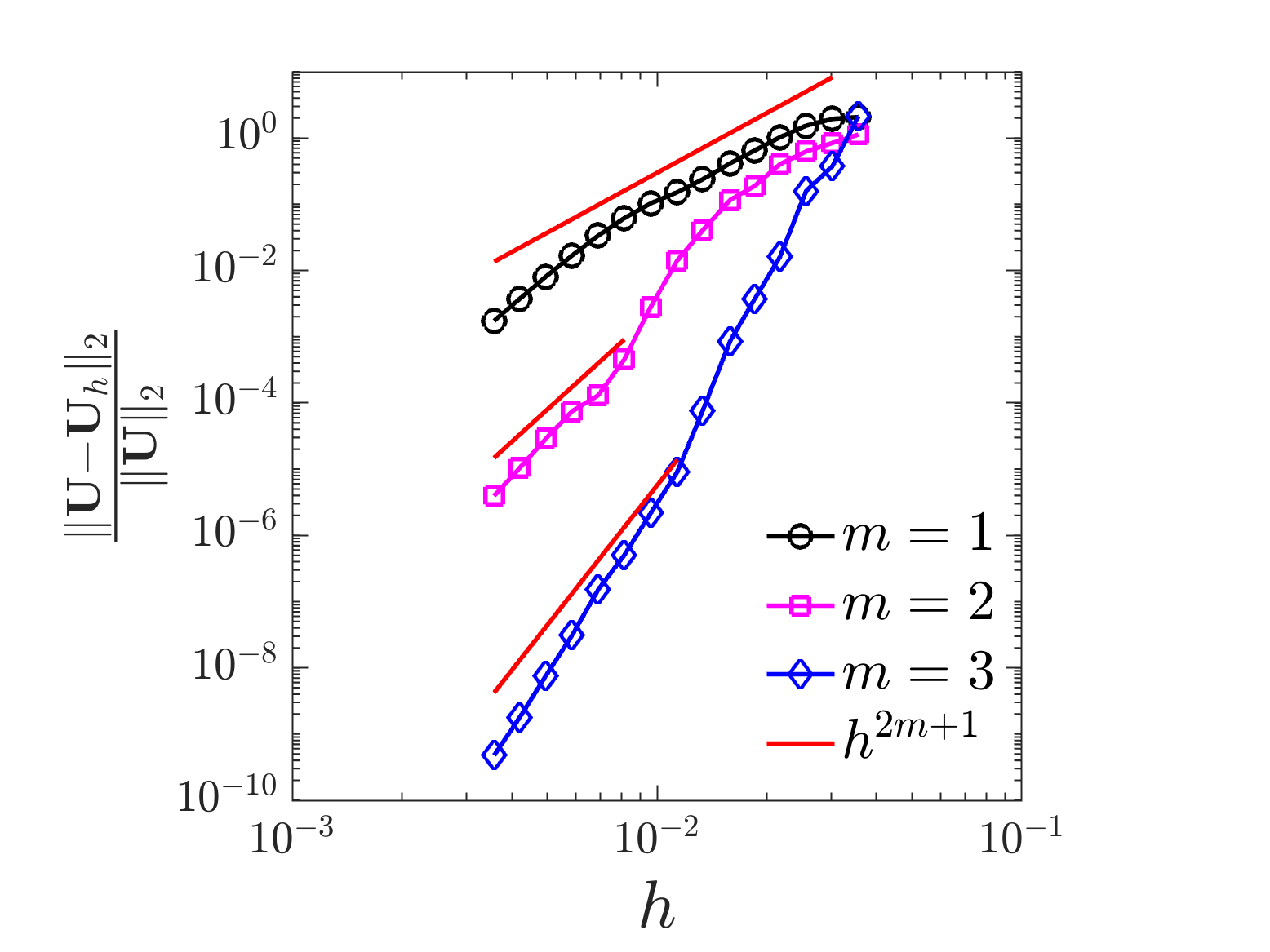}
    \includegraphics[width=\linewidth,trim={0cm 0cm 1.75cm 0cm},clip]{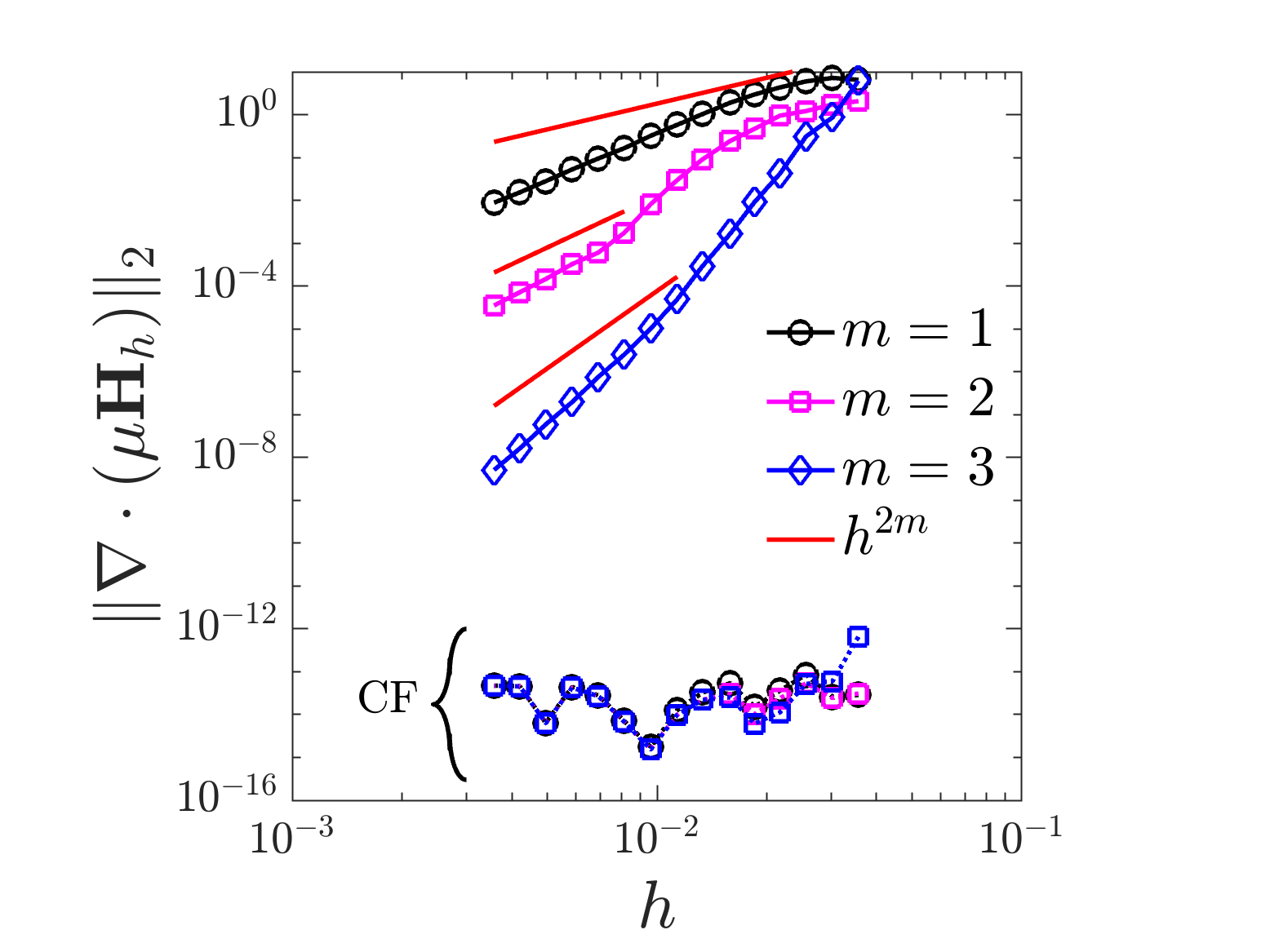}
        
    \end{adjustbox}
    \caption{PEC}
    \label{fig:conv_pec}
\end{subfigure}
\hfill
\begin{subfigure}{1.0\textwidth}
    \centering
    \begin{adjustbox}{max width=0.9\textwidth,center}
    \includegraphics[width=\linewidth,trim={0cm 0cm 1.75cm 0cm},clip]{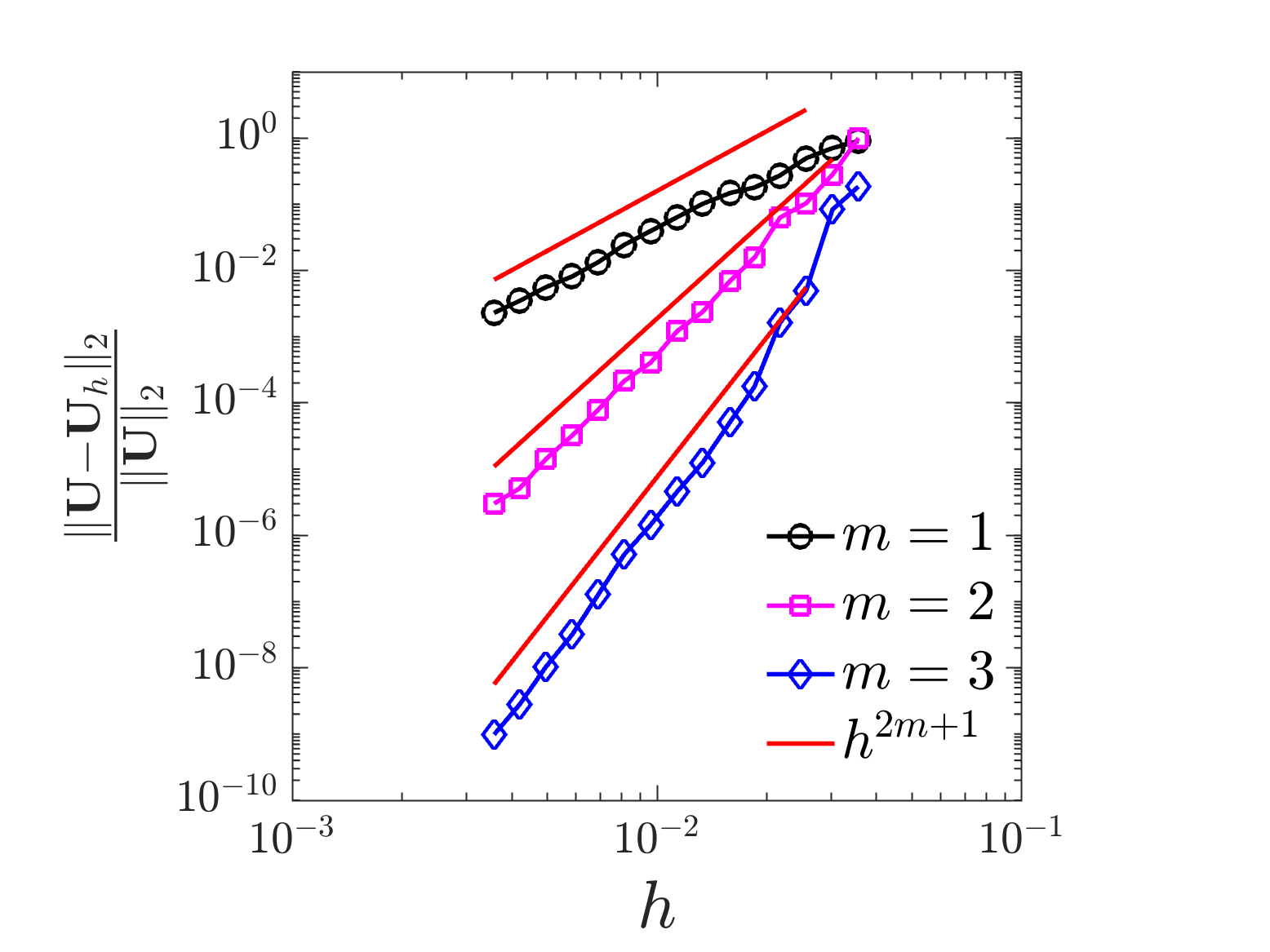}
    \includegraphics[width=\linewidth,trim={0cm 0cm 1.75cm 0cm},clip]{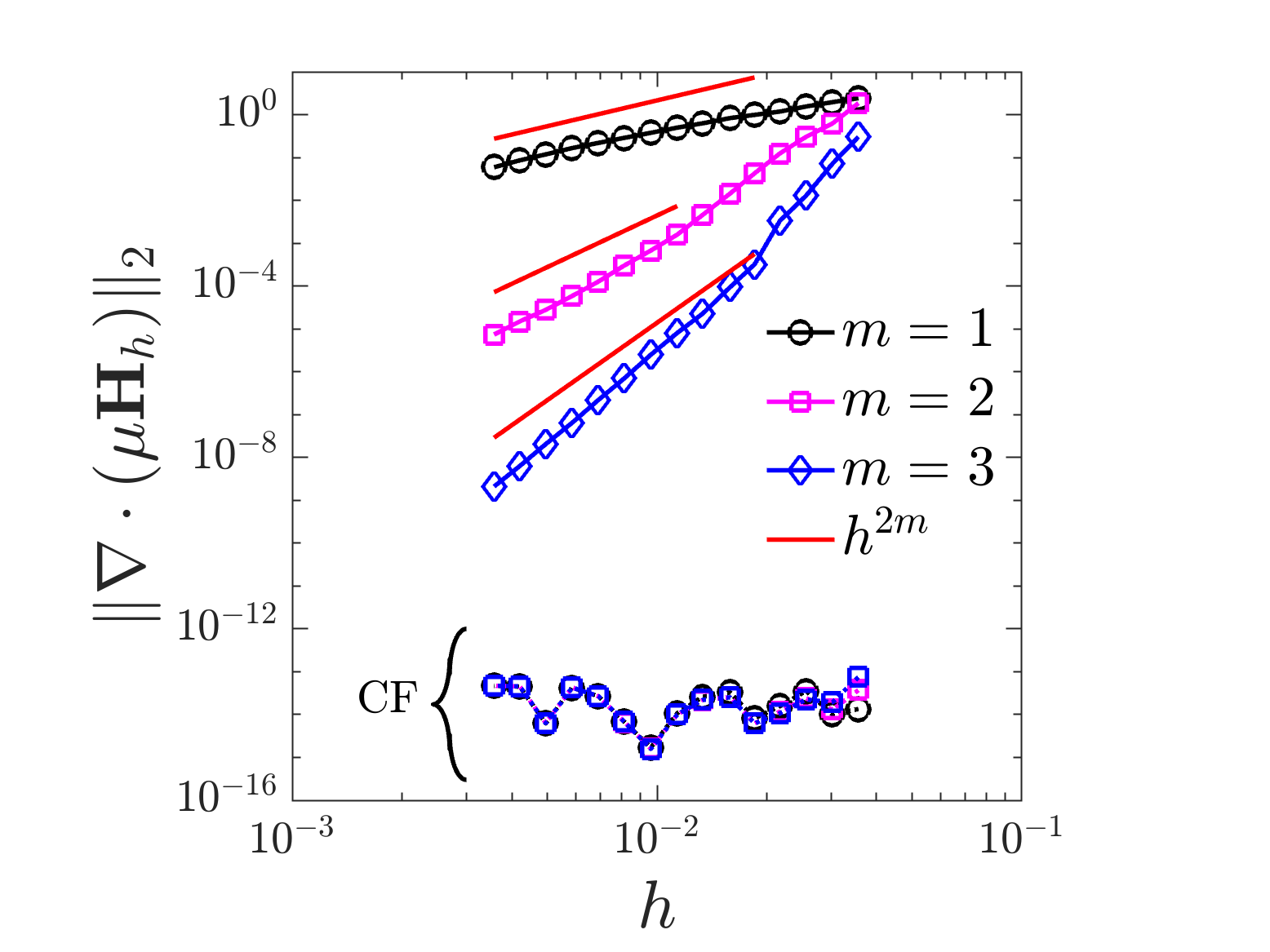}
        
    \end{adjustbox}
    \caption{PMC}
    \label{fig:conv_pmc}
\end{subfigure}

\vspace{1em}

\begin{subfigure}{1.0\textwidth}
    \centering
    \begin{adjustbox}{max width=0.9\textwidth,center}
    \includegraphics[width=\linewidth,trim={0cm 0cm 1.75cm 0cm},clip]{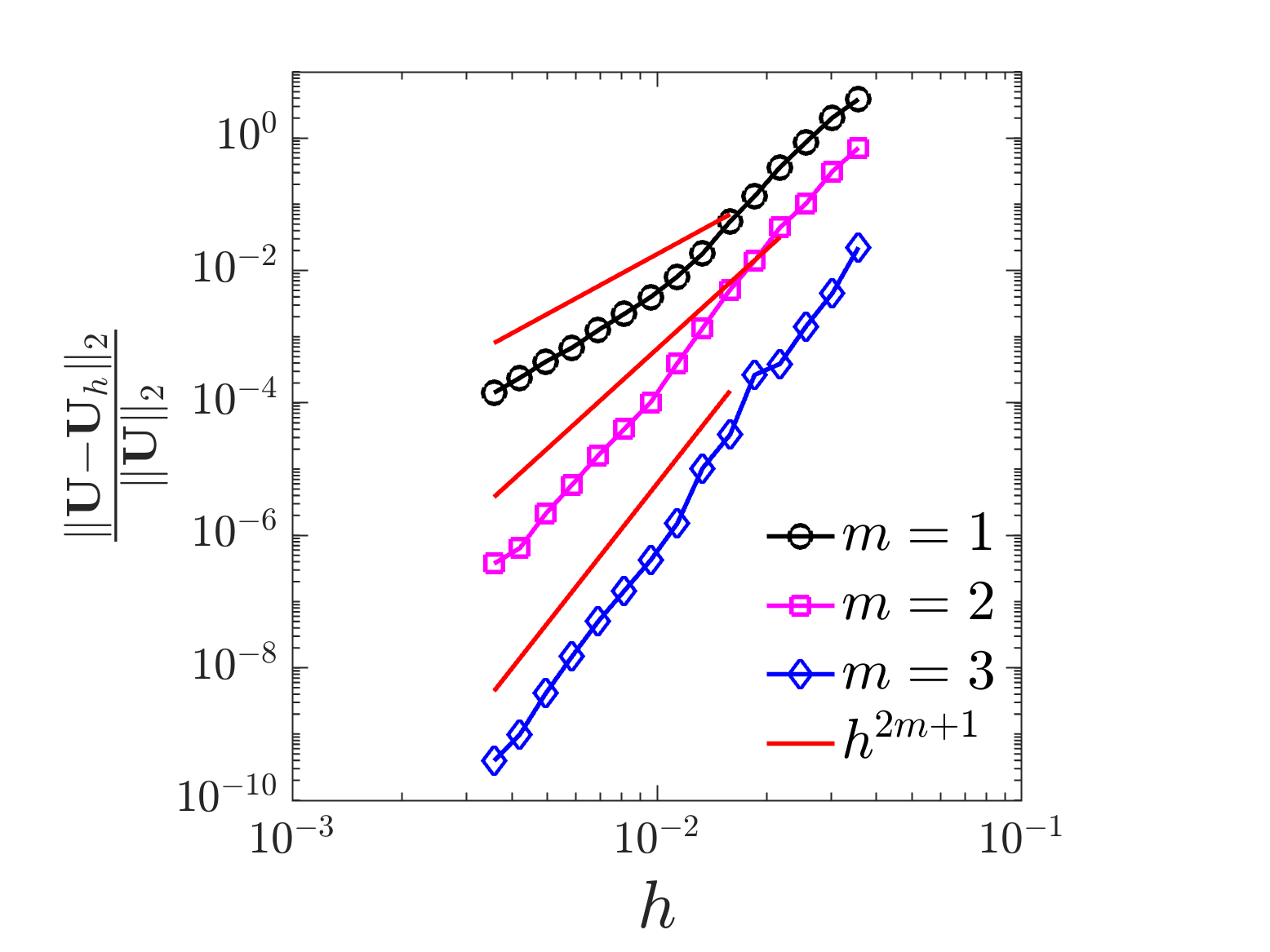}
    \includegraphics[width=\linewidth,trim={0cm 0cm 1.75cm 0cm},clip]{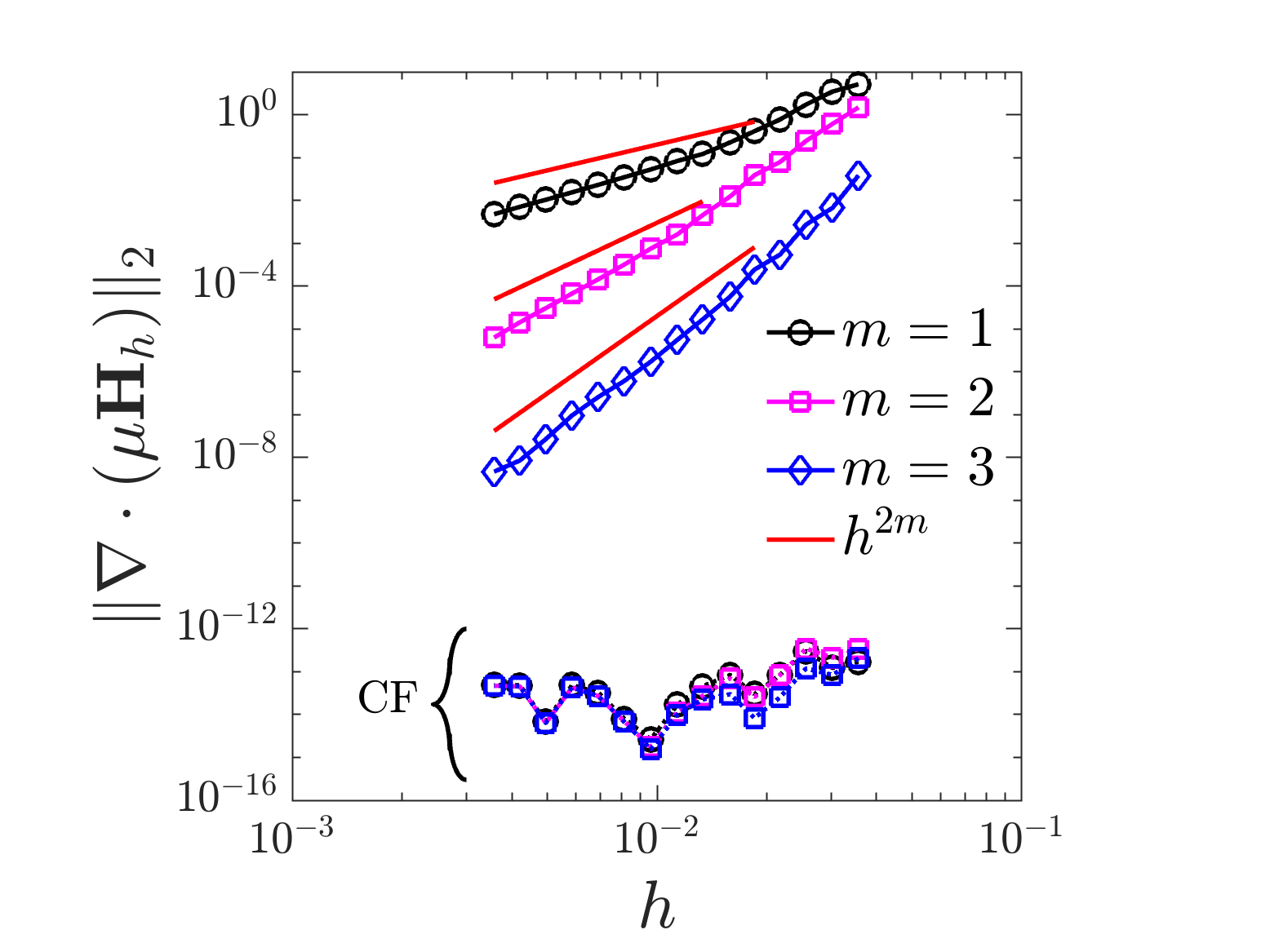}
        
    \end{adjustbox}
    \caption{Impedance}
    \label{fig:conv_impedance}
\end{subfigure}

\caption{Convergence plots for problems using PEC, PMC and impedance boundary conditions and different values of $m$. 
The left plots illustrate the relative error of $\mathbold{U} = [H_x; H_y; E_z]$ in 2-norm, 
    while the right plots show the norm \eqref{eq:computation_divergence_H} for the divergence of the magnetic field, 
    as well as the contribution of correction functions.}
\end{figure}

Convergence studies are illustrated in   Fig.~\ref{fig:conv_pec}, 
    Fig.~\ref{fig:conv_pmc} and Fig.~\ref{fig:conv_impedance} 
    for respectively PEC, 
    PMC and impedance boundary conditions. 
The time interval is $I = [0,1]$ and the mesh size varies between $1/280$ and $1/28$.
In all cases, 
    we obtain at least a $2m+1$ convergence order in 2-norm for all values of $m$ and all types of boundary conditions. 
Regarding the $L_2$-norm of the divergence of the magnetic field, 
    the right plots of Fig.~\ref{fig:conv_pec}, 
    Fig.~\ref{fig:conv_pmc} and Fig.~\ref{fig:conv_impedance} illustrate the norm \eqref{eq:computation_divergence_H},
    as well as the contribution of correction functions, 
    as a function of $h$.
The error in the divergence of correction functions approximating the magnetic field is due to the use of floating-point arithmetic, 
which is expected since we choose divergence-free correction function polynomials.

\begin{figure}   
	\centering
	\begin{adjustbox}{max width=0.9\textwidth,center}
		\includegraphics[width=\linewidth,trim={0cm 0cm 1.75cm 0cm},clip]{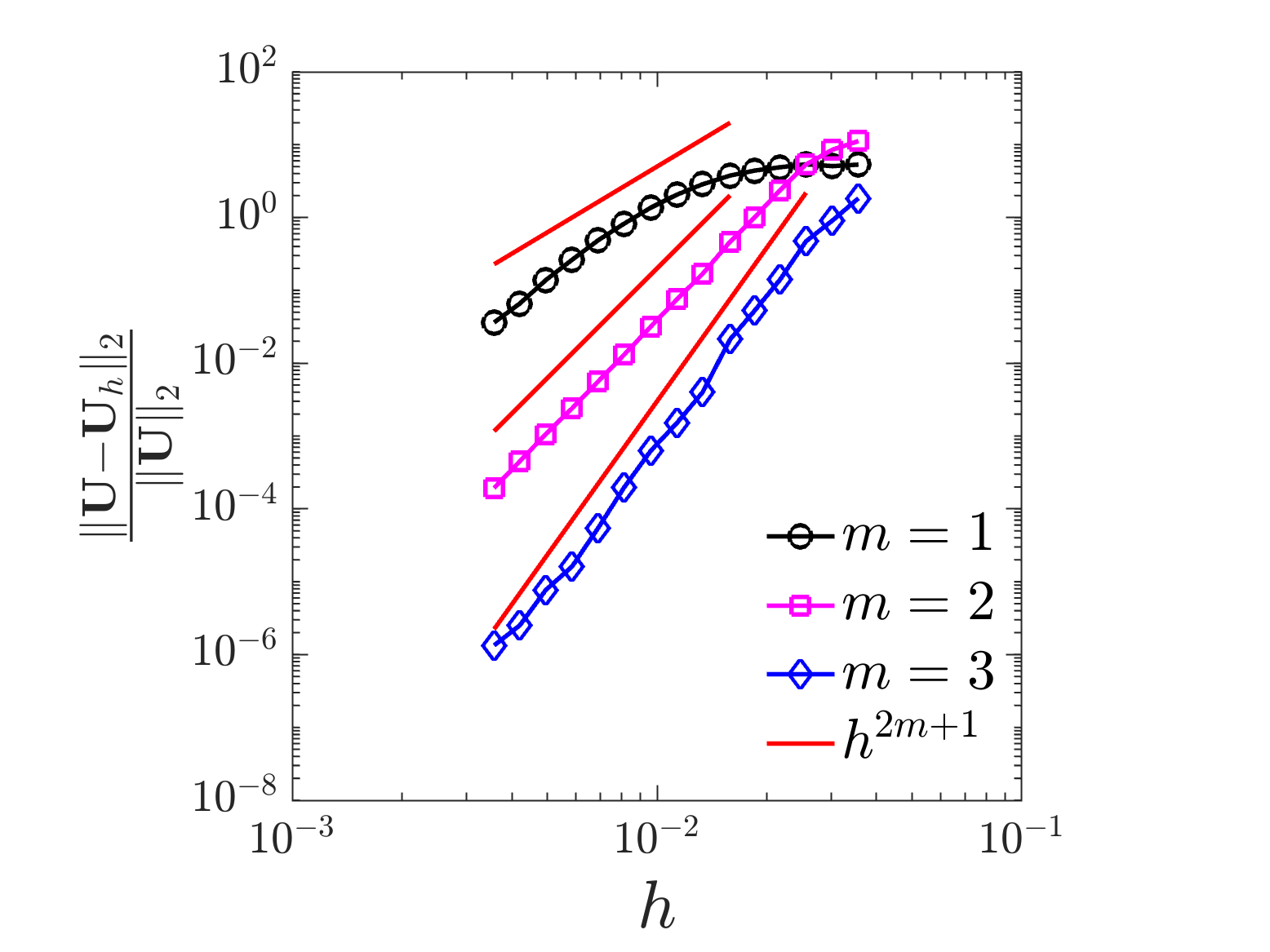} \hspace{-18.0pt}
		\includegraphics[width=\linewidth,trim={0cm 0cm 1.75cm 0cm},clip]{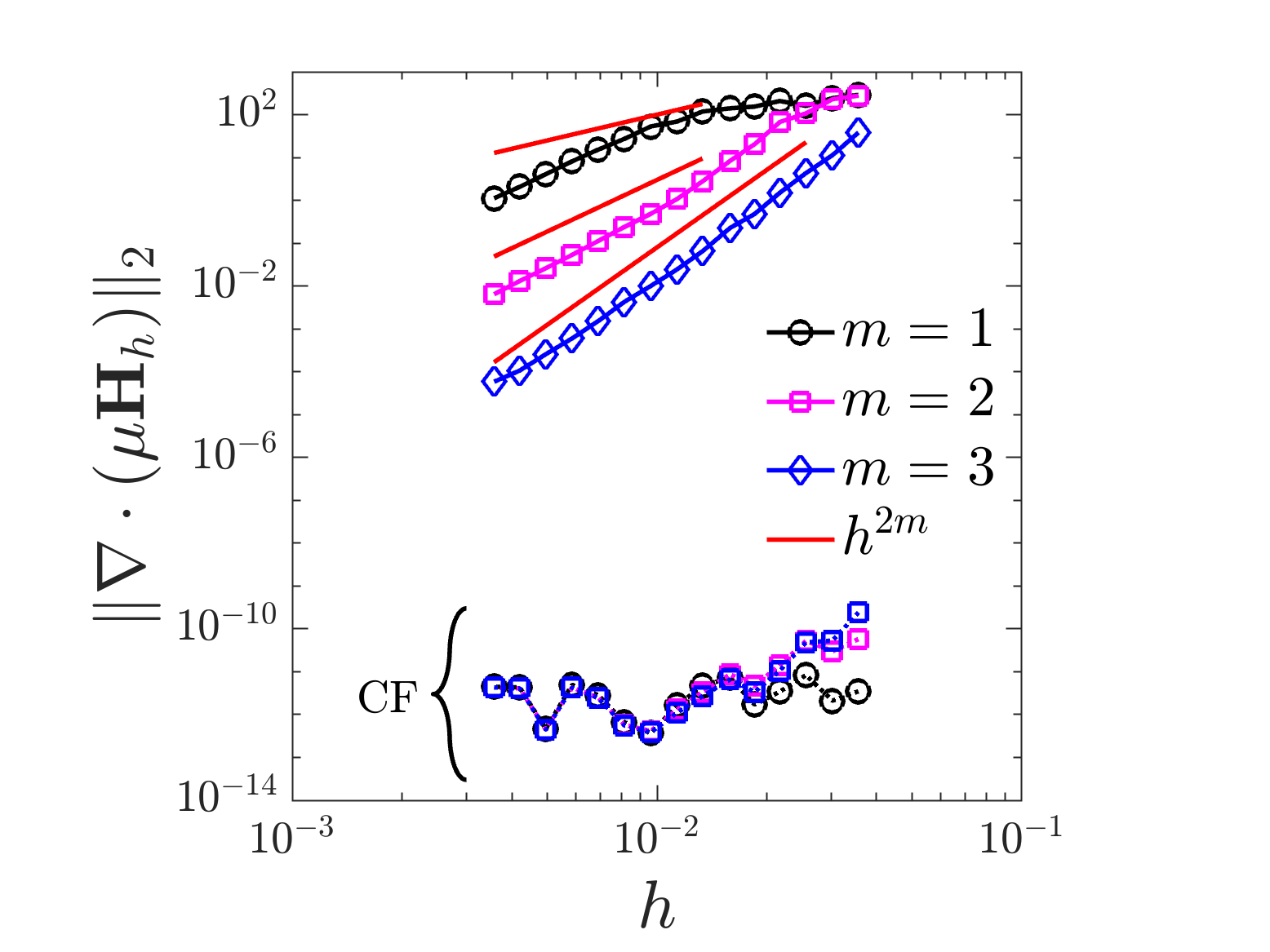}\hspace{-18pt}
	\end{adjustbox} 
       \caption{Convergence plots for a problem with variable coefficients using PEC condition and different values of $m$. 
The left plots illustrate the relative error of $\mathbold{U} = [H_x; H_y; E_z]$ in 2-norm, 
    while the right plots show the norm \eqref{eq:computation_divergence_H} for the divergence of the magnetic field, 
    as well as the contribution of correction functions.}
\label{fig:conv_pec_var_coeff}
\end{figure}

Let us consider the physical parameters to be dependent on the spatial coordinates $(x,y)$ that are given by 
\begin{equation}
    \mu(x,y) = e^{x+y}, \qquad \epsilon(x,y) = \frac{2(\omega^2+1)}{\omega^2}e^{-x-y}. 
\end{equation}
    with $\omega = 15\pi$.
The wave speed is then $c = \omega/\sqrt{2(\omega^2+1)} \approx 0.71$.
We choose the initial condition and PEC boundary condition so that the solution is given by 
\begin{equation}
    \begin{aligned}
    & H_x = \sin(\omega(x+y+t)), \\ 
    & H_y = -\sin(\omega(x+y+t)), \\
    & E_z = -\frac{\omega e^{x+y}}{\omega^2+1} \big( \omega\sin(\omega(x+y+t)) + \cos(\omega(x+y+t))\big).
    \end{aligned}
\end{equation}
The physical and computational domains are the same as before. 
Note that source terms are required for Amp\`ere-Maxwell's law \eqref{eq:Ampere_Maxwell_law} and PEC boundary condition \eqref{eq:2D_PEC_bnd_cdn}. 
We set $d=2m$, 
    $N_d = 2$ and the CFL constant to $0.5$ for all values of $m$.
Fig.~\ref{fig:conv_pec_var_coeff} shows a $2m+1$ order of convergence for the relative error on the electromagnetic fields. 
The expected $2m$ order of accuracy is also observed for the divergence-free constraint $\nabla\cdot(\mu\mathbold{H}) = 0$. 
Correction functions for the magnetic field are also divergence-free,  
    up to round-off errors. 

\begin{figure}   
	\centering
	\begin{adjustbox}{max width=1.0\textwidth,center}
	   \includegraphics[width=2.5in,trim={2.5cm 0cm 1.75cm 0cm},clip]{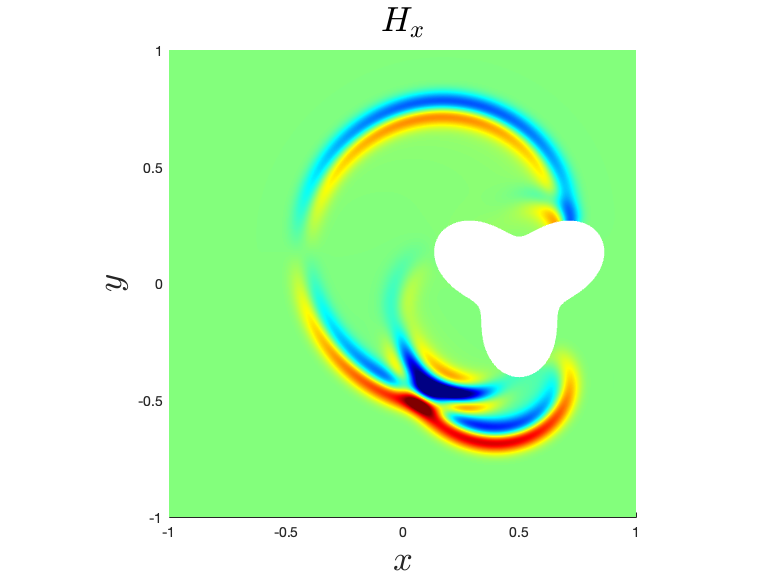} \hspace{-12.0pt}
		\includegraphics[width=2.5in,trim={2.5cm 0cm 1.75cm 0cm},clip]{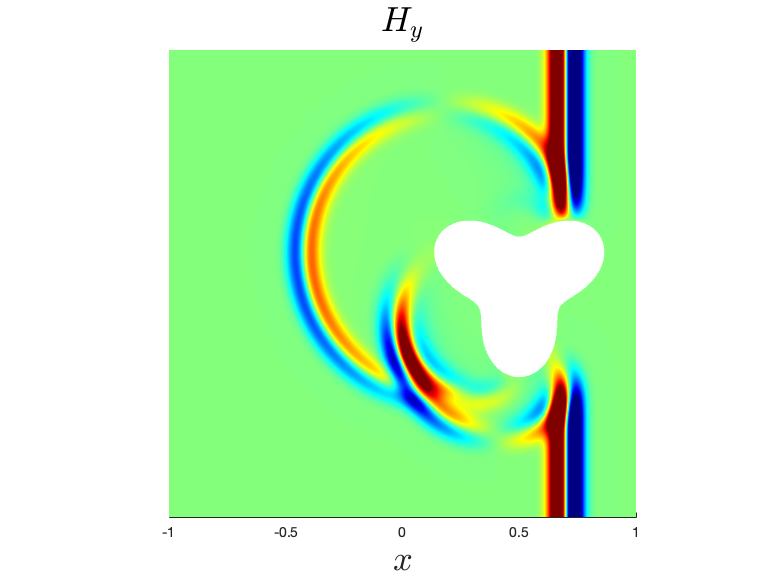} 
		\includegraphics[width=2.5in,trim={2.5cm 0cm 1.75cm 0cm},clip]{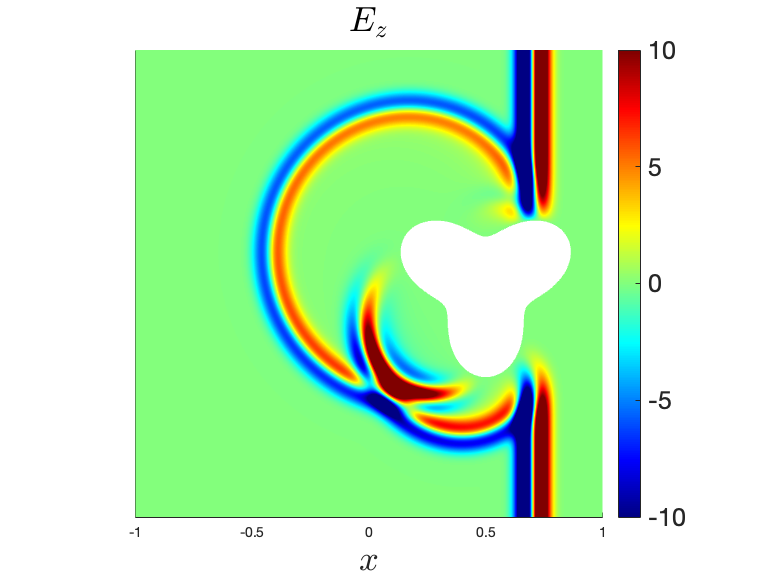} 
	\end{adjustbox} 
       \caption{Reference solution $\mathbold{U}^*$ at the final time for self-convergence studies in the PEC case.}
\label{fig:self_conv_pec_ref_sol}
\end{figure}

\begin{figure}   
	\centering
	\begin{adjustbox}{max width=0.9\textwidth,center}
		\includegraphics[width=\linewidth,trim={0cm 0cm 1.75cm 0cm},clip]{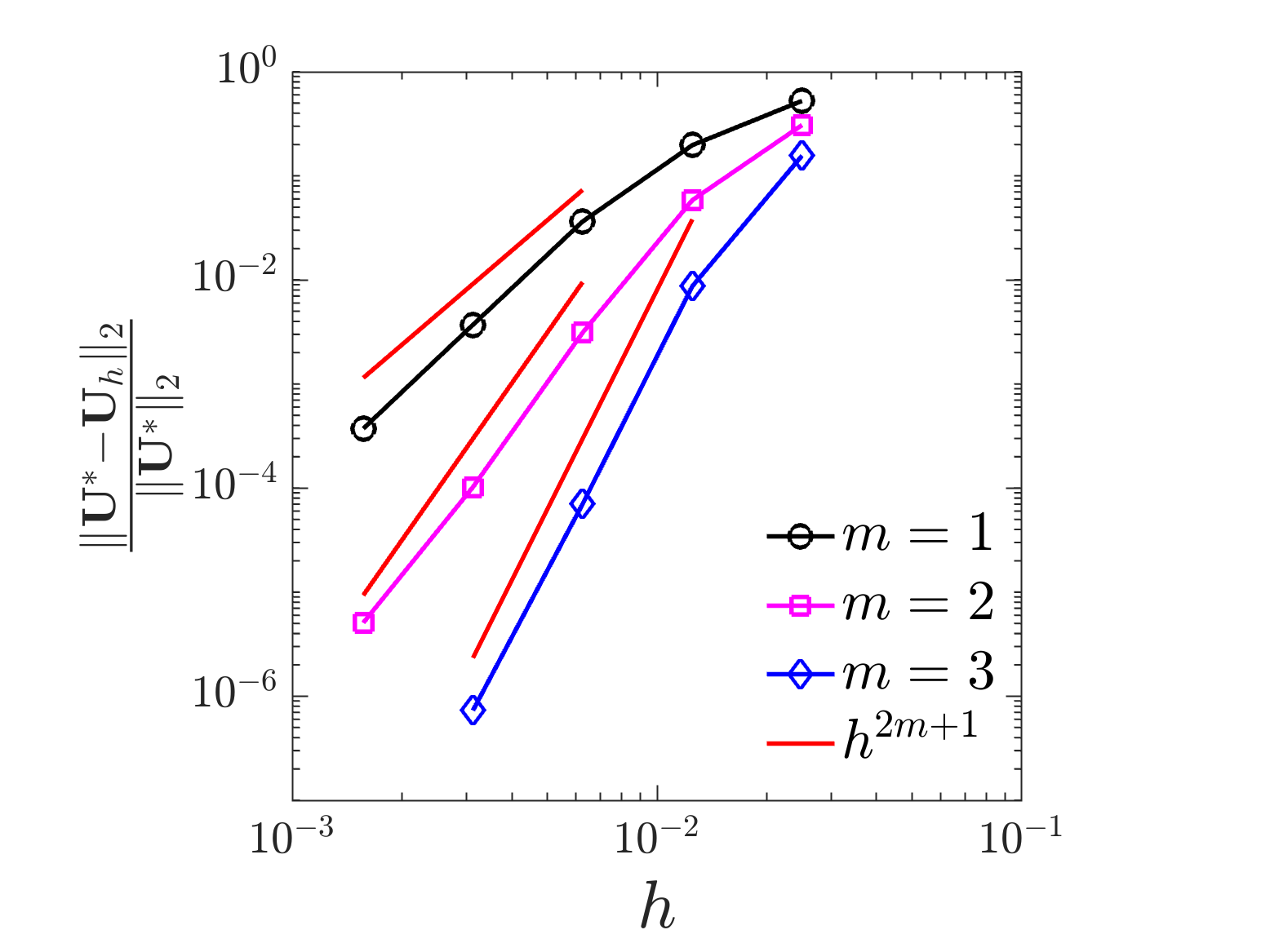} \hspace{-18.0pt}
		\includegraphics[width=\linewidth,trim={0cm 0cm 1.75cm 0cm},clip]{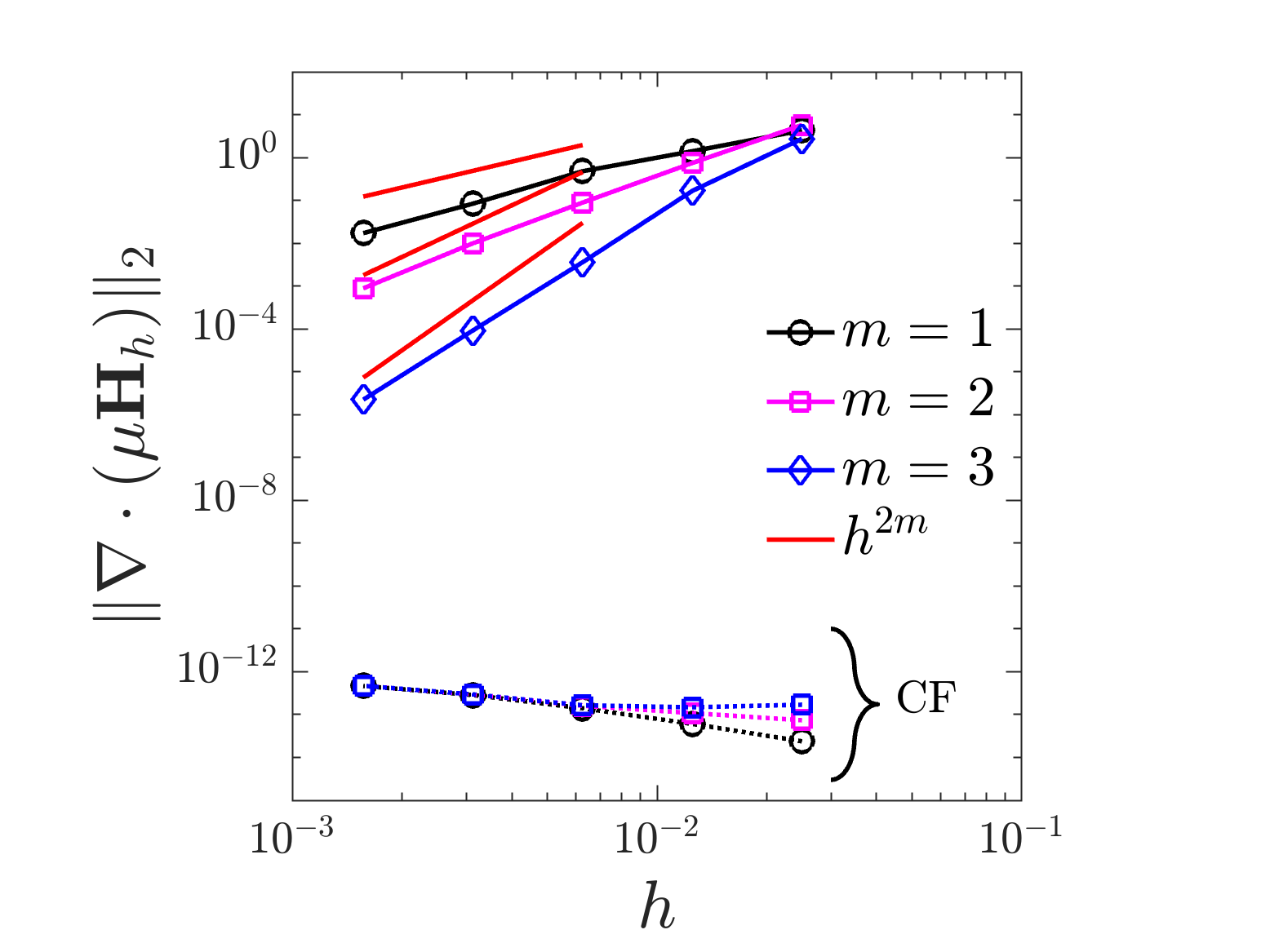}\hspace{-18pt}
	\end{adjustbox} 
       \caption{Self-convergence plots for the illumination of a PEC by the incident wave \eqref{eq:incident_wave} and different values of $m$. 
The left plots illustrate the relative error of $\mathbold{U} = [H_x; H_y; E_z]$ in 2-norm using the reference solution $\mathbold{U}^*$, 
    while the right plots show the norm \eqref{eq:computation_divergence_H} for the divergence of the magnetic field, 
    as well as the contribution of correction functions.}
\label{fig:self_conv_pec}
\end{figure}

We now consider a PEC 
    illuminated by the incident wave 
\begin{equation} \label{eq:incident_wave}
H_x =0, \qquad H_y = -\frac{2}{\sigma^2}(x-\gamma-t)e^{-\frac{(x-\gamma-t)^2}{\sigma^2}}, \qquad E_z = -H_y,
\end{equation}
    with $\gamma = -0.3$ and $\sigma=1/20$.
The PEC has a smooth three-star shape 
    described by \eqref{eq:level_set_fct} with $x_0 =0.5$ and $y_0 =0$.
The computational domain is $\Omega = [-1,1]\times[-1,1]$ with periodic boundary conditions and $I = [0,1]$.
The physical parameters are $\mu =1$ and $\epsilon =1$.
To the best of our knowledge,
    there is no known analytical solution to this problem. 
We set the CFL constant to be $0.5$ for all values of $m$.
We therefore perform self-convergence studies with a reference solution obtained with $m=3$ and $h=1/1280$.
Fig.~\ref{fig:self_conv_pec_ref_sol} illustrates the reference solution at the final time. 
The coarser meshes with $h\in\{1/80,1/160,1/320,1/640\}$ are chosen such that all their primal nodes are  contained within the mesh of the reference solution.
Fig.~\ref{fig:self_conv_pec} illustrates the self-convergence plot for the electromagnetic fields and the error on the divergence of the magnetic field. 
In all cases, 
    we observe the expected order of convergence.

\subsection{Interface problems}

Consider now the domain $\Omega$ subdivided into two subdomains $\Omega^+$ and $\Omega^-$, 
    where $\Omega^-$ contains $\Omega^+$, which is delimited by the interface $\tilde{\Gamma}$, 
    defined by the level set function \eqref{eq:level_set_fct}.
Periodic boundary conditions are enforced on the boundary of $\Omega^-$.
The physical parameters are $\mu^- = 3$,
    $\epsilon^- = 2/3$, 
    $\mu^+= 1$ and $\epsilon^+ = 1$.
The initial, 
    boundary and interface conditions are chosen so that the solution is 
\begin{equation}
    H_x^- = H_y^- = -\cos(2\pi\alpha(x+y+t))/3, \qquad
    E_z^- = \cos(2\pi\alpha(x+y+t)),
\end{equation}
    with $\alpha =3$ in $\Omega^-$, 
    and $H^+$,
    $H_y^+$ and $E_z^+$ are given by \eqref{eq:sol_bnd_cdn_cst_coeff} in $\Omega^+$. 
Note that source terms are required for the interface conditions \eqref{eq:interface_cdns_2d}.

\begin{figure}   
	\centering
	\begin{adjustbox}{max width=1.0\textwidth,center}
		\includegraphics[width=2.5in,trim={0cm 0cm 1.75cm 0cm},clip]{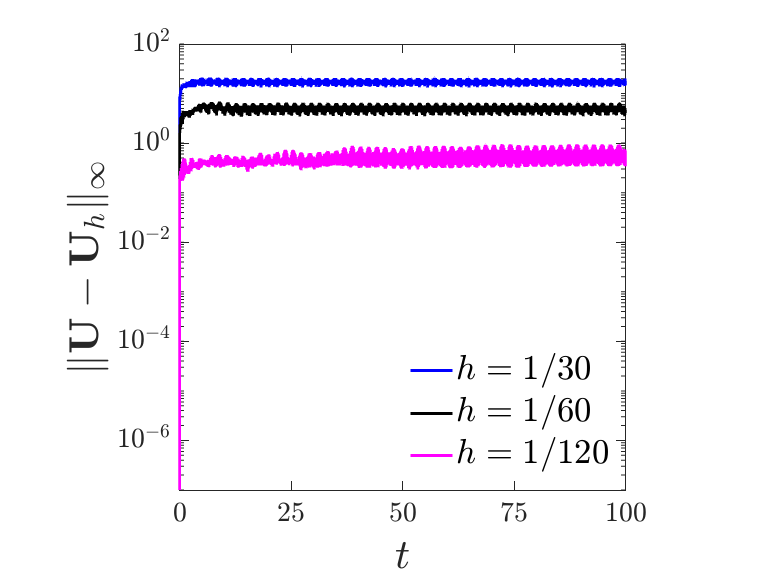} \hspace{-18.0pt}
		\includegraphics[width=2.5in,trim={0cm 0cm 1.75cm 0cm},clip]{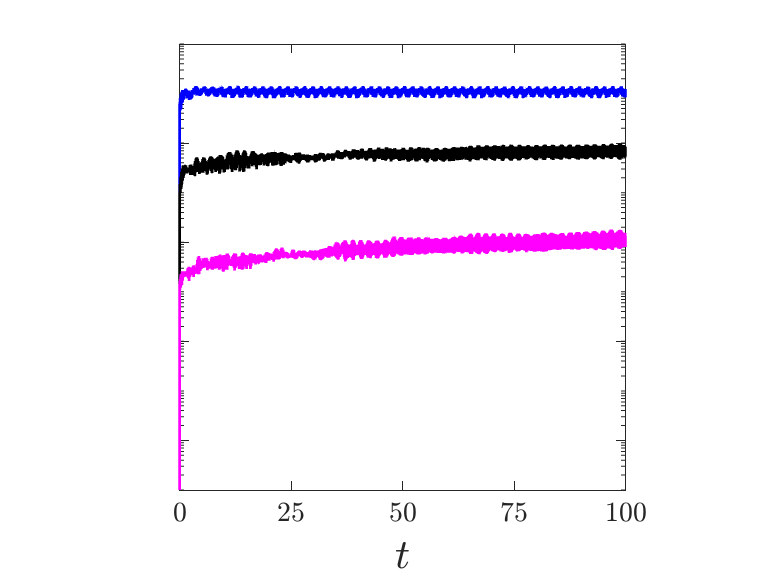}\hspace{-18pt}
		\includegraphics[width=2.5in,trim={0cm 0cm 1.75cm 0cm},clip]{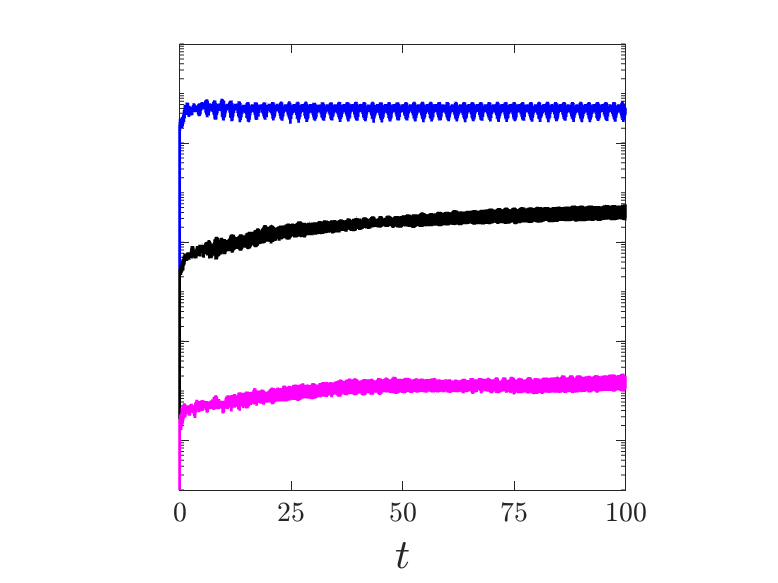}\hspace{-18pt}
	\end{adjustbox} 
       \caption{Maximum norm of the error as a function of time for a problem using interface conditions \eqref{eq:interface_cdns_2d} with source terms, different mesh sizes and values of $m$. The columns are for different $m$: 1 to 3 from left to right.
       Here $\mathbold{U} = [H_x;H_y;E_z]$.}
\label{fig:long_time_simulations_interface_cdn}
\end{figure}

\begin{figure}   
	\centering
	\begin{adjustbox}{max width=0.9\textwidth,center}
		\includegraphics[width=\linewidth,trim={0cm 0cm 1.75cm 0cm},clip]{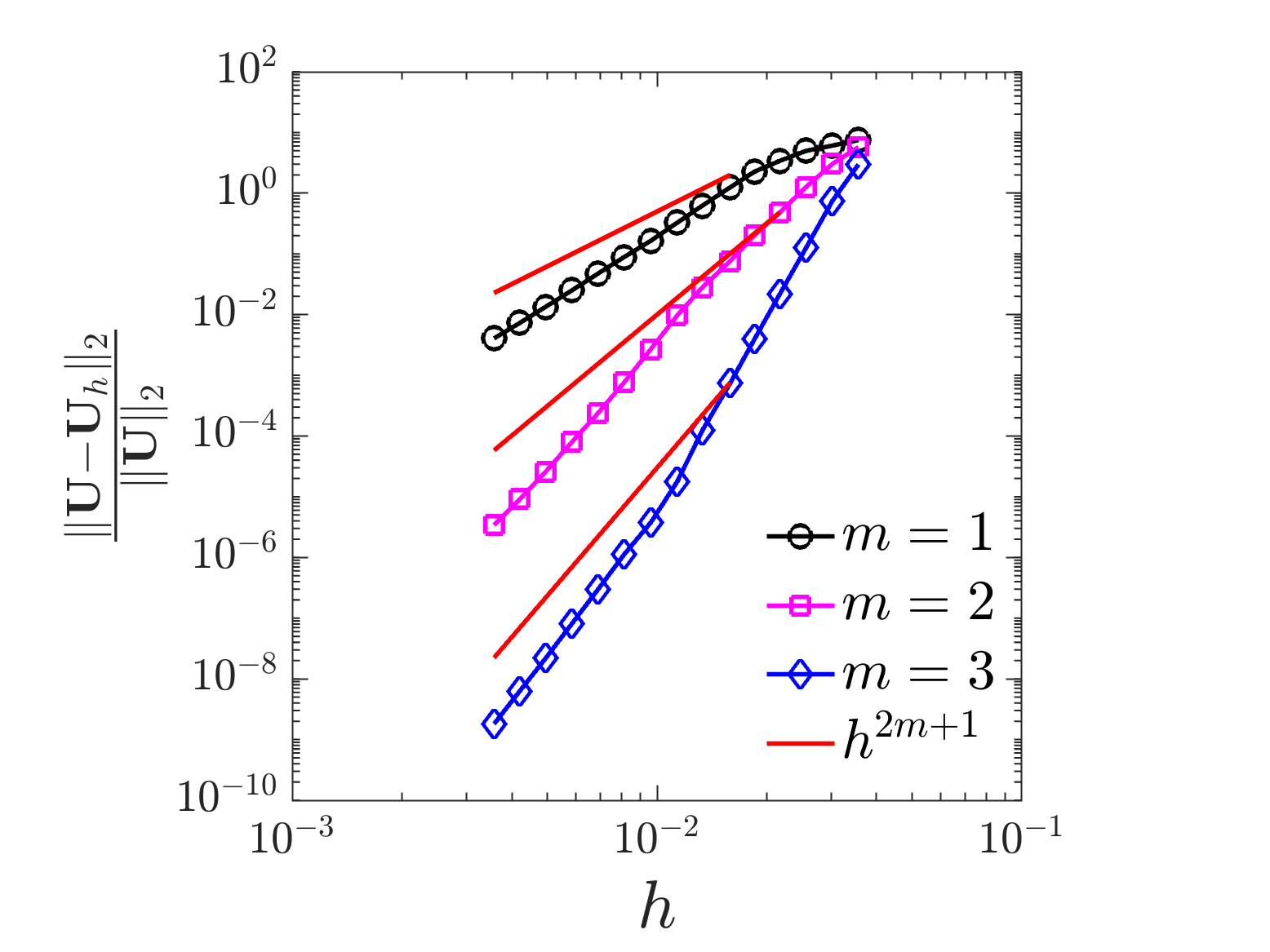} \hspace{-18.0pt}
		\includegraphics[width=\linewidth,trim={0cm 0cm 1.75cm 0cm},clip]{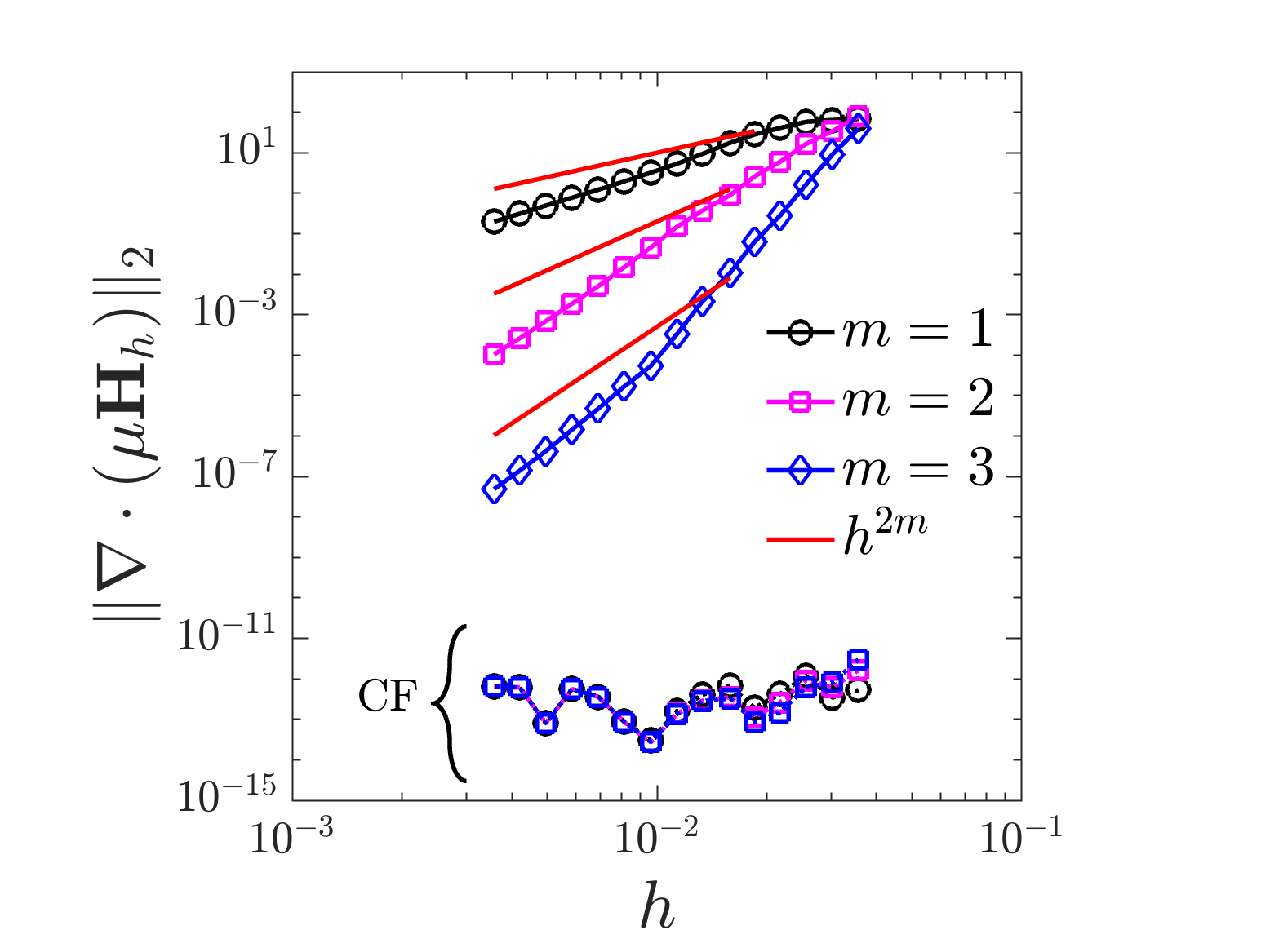}\hspace{-18pt}
	\end{adjustbox} 
	\begin{adjustbox}{max width=0.9\textwidth,center}
		\includegraphics[width=\linewidth,trim={0cm 0cm 1.75cm 0cm},clip]{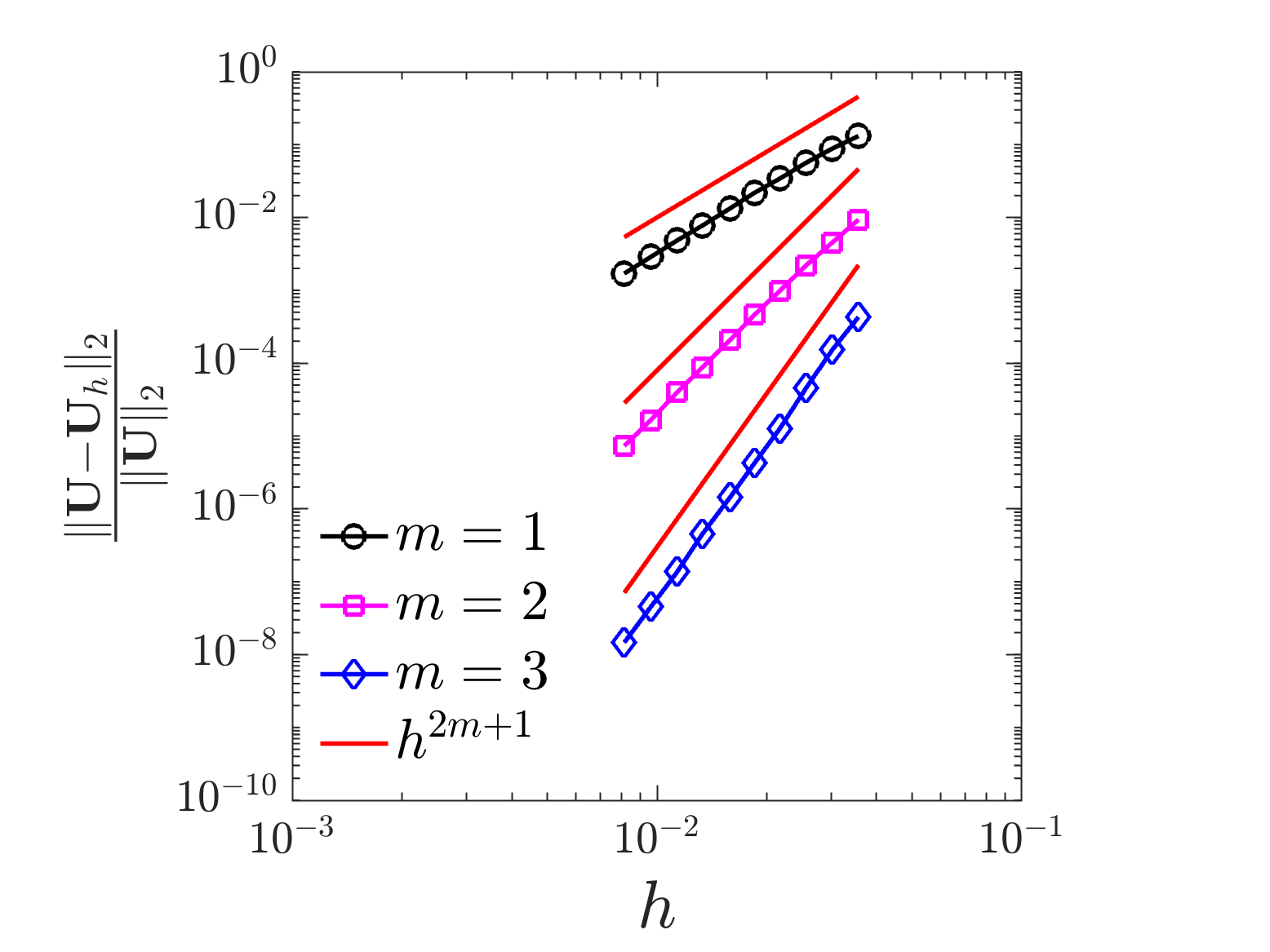} \hspace{-18.0pt}
		\includegraphics[width=\linewidth,trim={0cm 0cm 1.75cm 0cm},clip]{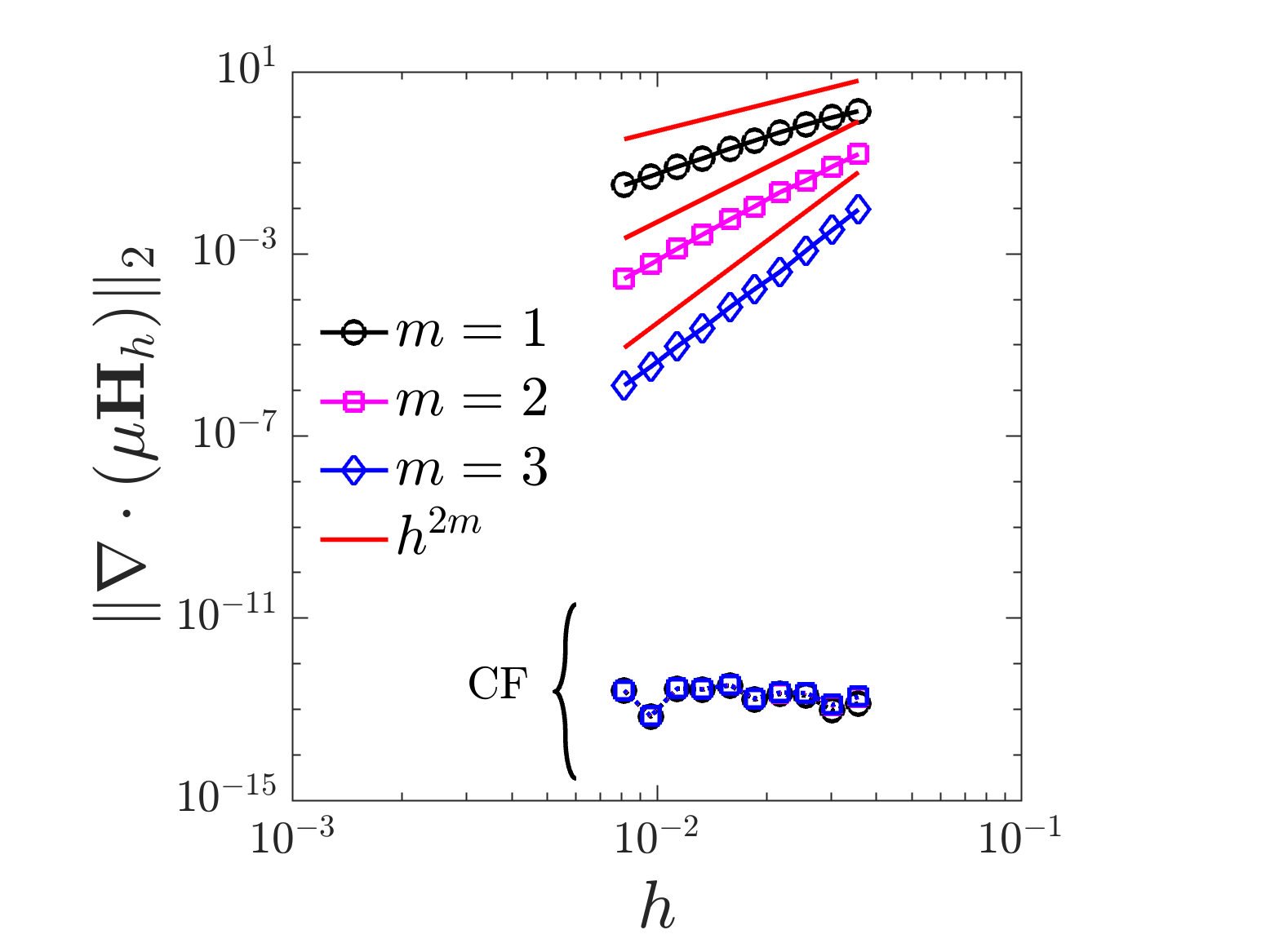}\hspace{-18pt}
	\end{adjustbox} 
       \caption{Convergence plots for problems using interface conditions and different values of $m$. 
The left plots illustrate the relative error of $\mathbold{U} = [H_x; H_y; E_z]$ in 2-norm, 
    while the right plots show the norm \eqref{eq:computation_divergence_H} for the divergence of the magnetic field, 
    as well as the contribution of correction functions.
The top row corresponds to a problem using a smooth three-star magnetic dielectric, 
while the bottom row is for a circular magnetic dielectric.}
\label{fig:conv_dielectric}
\end{figure}

As for the boundary condition case,    
    we perform long-time simulations to investigate the stability of the Hermite-Taylor discrete correction function method for interface problems. 
We consider $I=[0,100]$,
    $h\in\{1/30,1/60,1/120\}$,
    $m=1-3$ and the CFL constant to be $0.5$ for all values of $m$. 
The parameters for the discrete correction function method are $d=2m$ and $N_d = 2m$.
Fig.~\ref{fig:long_time_simulations_interface_cdn} shows the error in maximum norm as a function of time. 
The results suggest that the method is stable.
The top row of Fig.~\ref{fig:conv_dielectric} illustrates the convergence plots for the electromagnetic fields and the divergence-free constraint on the magnetic field at the final time $t_f =1$. 
The errors converge at least to the expected order, 
    that is, 
    $2m+1$ and $2m$ for the electromagnetic fields using the 2-norm and the divergence of the magnetic field using the norm \eqref{eq:computation_divergence_H}, 
    respectively.
The maximum condition number of DCFM matrices for all three meshes is approximately $10^5$,
    $10^7$ and $10^9$ for $m=1-3$, 
    respectively.

As a second interface problem, 
    we consider a magnetic dielectric cylinder in free space illuminated by a time-harmonic incident plane wave. 
The exact solution to this problem can be found in \cite{Umashankar1993,Cai2003}.
The computational domain is $\Omega = [-1,1]\times[-1,1]$, 
    and is divided into two subdomains: 
    $\Omega^+$ and $\Omega^-$, representing the free space and the dielectric, 
    respectively.
The interface between the subdomains is represented by a circle centered at $(0,0)$ with a radius of $0.6$, 
    which encloses the subdomain $\Omega^-$.
We also consider a boundary of $\Omega^+$ consisting of a circle of radius $0.8$, 
    also centered at $(0,0)$, 
    on which PEC boundary condition \eqref{eq:2D_PEC_bnd_cdn} with source terms is enforced. 
The time interval is $I = [0,1]$.
The physical parameters are $\mu^+ = 1$, 
    $\epsilon^+=1$, 
    $\mu^- = 2$ and $\epsilon^- = 2.25$.
The parameters of the DCFM are the same as those of the previous problem.
The bottom row of Fig.~\ref{fig:conv_dielectric} illustrates the convergence plots for the electromagnetic fields and the divergence-free constraint on the magnetic field. 
The results are in agreement with the theory. 

\begin{figure}   
	\centering
	\begin{adjustbox}{max width=1.0\textwidth,center}
	   \includegraphics[width=2.5in,trim={2.5cm 0cm 1.75cm 0cm},clip]{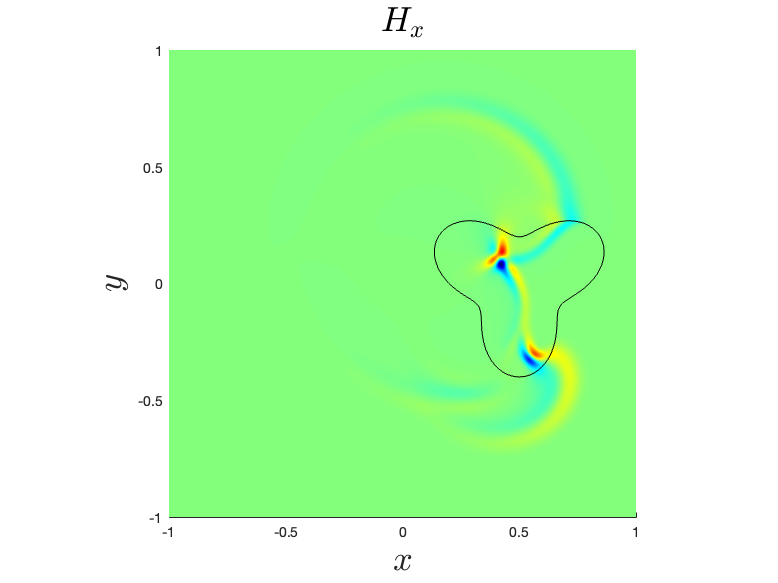} \hspace{-12.0pt}
		\includegraphics[width=2.5in,trim={2.5cm 0cm 1.75cm 0cm},clip]{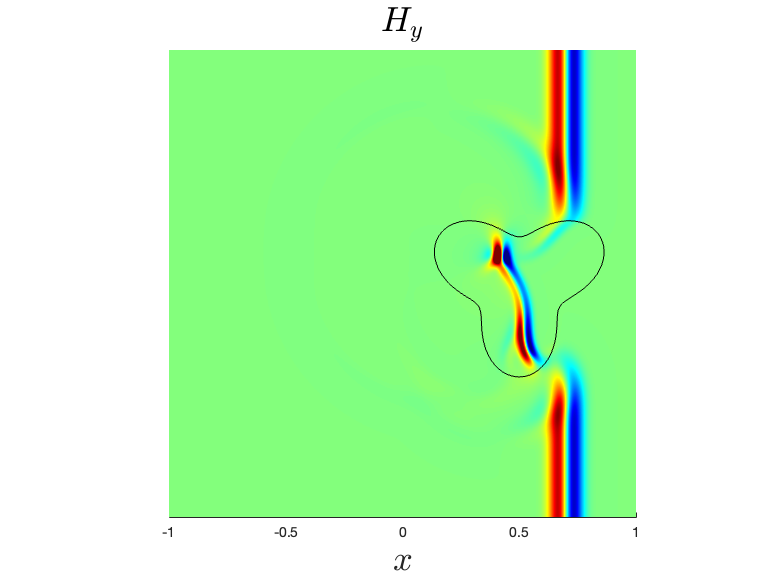} 
		\includegraphics[width=2.5in,trim={2.5cm 0cm 1.75cm 0cm},clip]{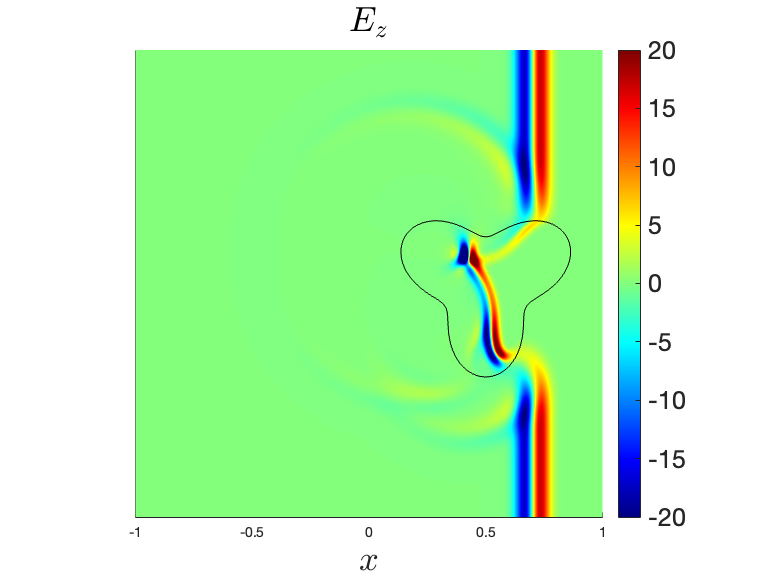} 
	\end{adjustbox} 
       \caption{Reference solution $\mathbold{U}^*$ at the final time for self-convergence studies in the dielectric case. The black line represents the interface $\tilde{\Gamma}$.}
\label{fig:self_conv_dielectric_ref_sol}
\end{figure}

\begin{figure}   
	\centering
	\begin{adjustbox}{max width=0.9\textwidth,center}
		\includegraphics[width=\linewidth,trim={0cm 0cm 1.75cm 0cm},clip]{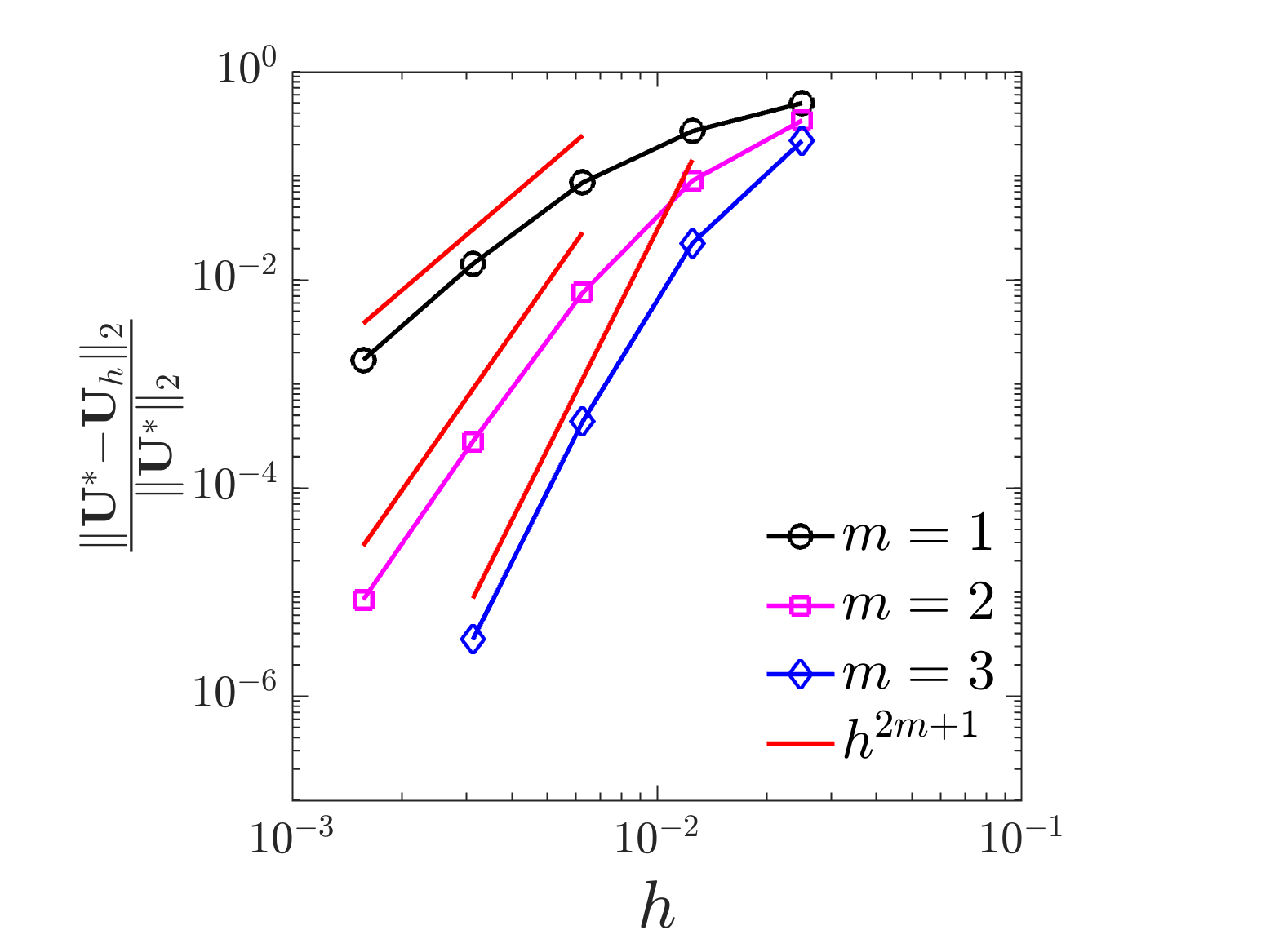} \hspace{-18.0pt}
		\includegraphics[width=\linewidth,trim={0cm 0cm 1.75cm 0cm},clip]{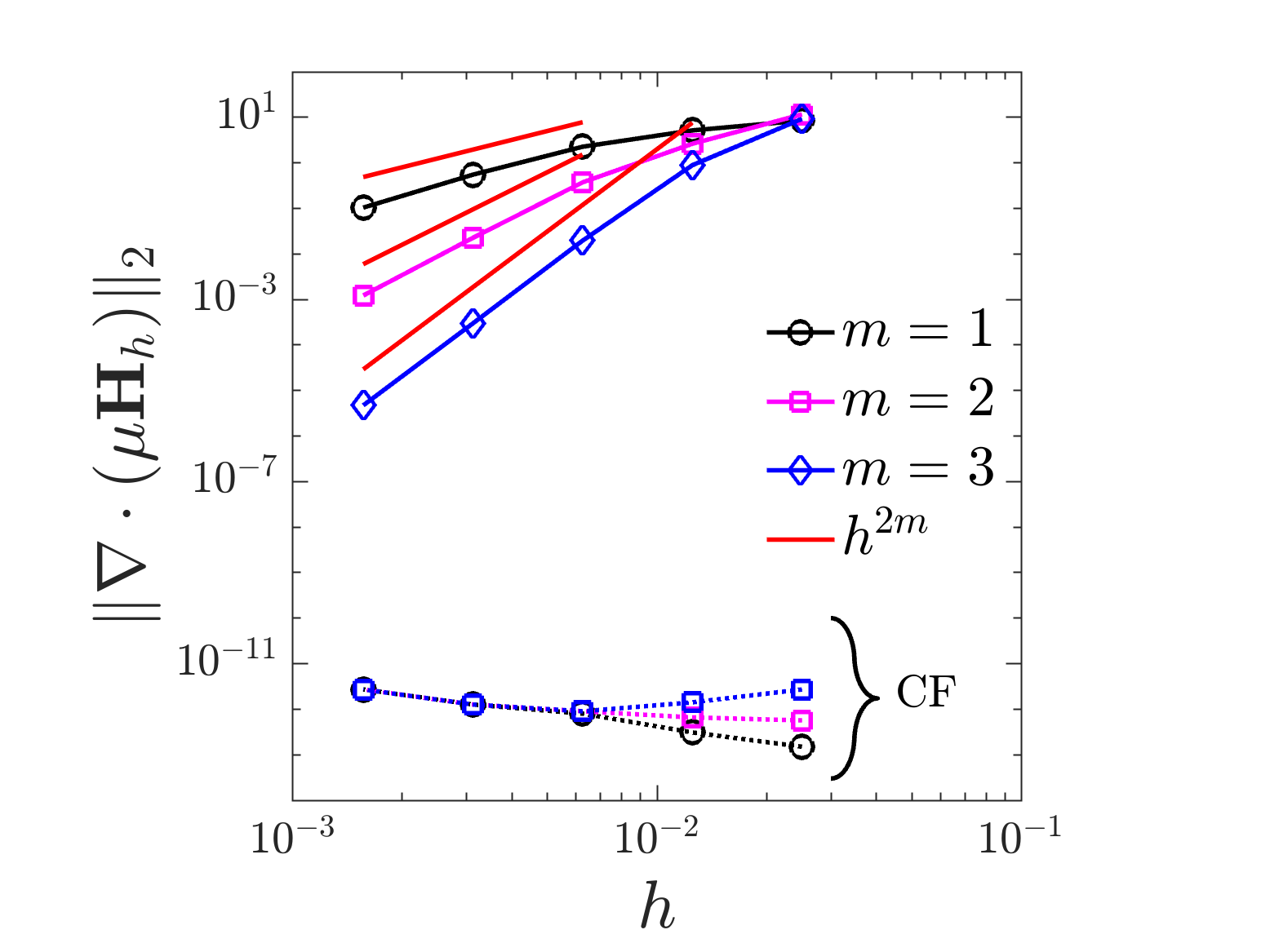}\hspace{-18pt}
	\end{adjustbox} 
       \caption{Self-convergence plots for the illumination of a magnetic dielectric by the incident wave \eqref{eq:incident_wave} and different values of $m$. 
The left plots illustrate the relative error of $\mathbold{U} = [H_x; H_y; E_z]$ in 2-norm using the reference solution $\mathbold{U}^*$, 
    while the right plots show the norm \eqref{eq:computation_divergence_H} for the divergence of the magnetic field, 
    as well as the contribution of correction functions.}
\label{fig:self_conv_dielectric}
\end{figure}

Let us now consider a dielectric with a smooth three-star shape illuminated by the incident wave \eqref{eq:incident_wave}.
The subdomain $\Omega^+$, 
    representing the dielectric,
    is contained in $\Omega^-$ and the interface $\tilde{\Gamma}$ is defined by the level set function \eqref{eq:level_set_fct}.
The parameters are the same as in the embedded boundary case, 
    except that the physical parameters are $\mu^- = 1$, 
    $\epsilon^-=1$, 
    $\mu^+ = 2$ and $\epsilon^+ = 2$.
The reference solution at the final time is illustrated in Fig.~\ref{fig:self_conv_dielectric_ref_sol}.
The self-convergence plot for the electromagnetic fields is illustrated in the left plot of Fig.~\ref{fig:self_conv_dielectric}, 
    while the right plot shows the error on the divergence of the magnetic field. 
The expected orders of convergence are observed
    and correction functions approximating $\mathbold{H}$ are divergence-free, 
    up to round-off errors. 
\subsection{Flat metasurface problems}

We now consider a flat vertical metasurface modeled by generalized sheet transition condition \eqref{eq:GSTC_flat} and located at $x_{\Gamma}=1/2+\pi/100$. 
The physical parameters of the metasurface are $\mu_0 =1$, 
    $\epsilon_0=1$, 
    $\chi_{mm}^{xx} = 2/(3\omega)$ and
    $\chi_{mm}^{yy} =\chi_{ee}^{zz} = 2/\omega$ with $\omega = 4 \pi$.
The properties of the medium surrounding either side of the metasurface are 
    $\mu=1$ and $\epsilon = 2$.
The wave speed is therefore $c=1/\sqrt{2}$.
The computational domain is $\Omega = [0,1]\times[0,1]$.
The subdomains are $\Omega^- = [0.2,x_{\Gamma})\times[0,1]$ and $\Omega^+ = (x_{\Gamma},0.8]\times[0,1]$.
We enforce periodic boundary conditions on the top and the bottom of the domain, 
    and PEC boundary conditions with source terms are enforced at $(0.2,y)$ and $(0.8,y)$ for all $y\in[0,1]$.
The initial condition and the source terms of PEC boundary condition are chosen so that the solution is 
\begin{equation}
    \begin{aligned}
    & H_x^- = -\sin(\omega(x+y+t)),\qquad H_y^- = E_z^- = \sin(\omega(x+y+t)), \\
    & H_x^+ = -\cos(\omega(x+y+t)),\qquad H_y^+ = E_z^+ =\cos(\omega(x+y+t)).
    \end{aligned}
\end{equation}
There is no source term for generalized sheet transition condition \eqref{eq:GSTC_flat}. 
If we set $d=2m$ and $N_d = 2m-1$,
    we expect a $2m+1$ order Hermite-Taylor discrete correction function method, 
    as in the previous numerical examples for interface and embedded boundary problems. 
However, 
    numerical results in this setting suggest that,
    although the method is stable according to long-time simulations,
    we observe only a $2m$ order of convergence. 
We therefore increase the degree of the polynomial correction functions
    to $d = 2m+1$ and set $N_d = 2m$ to obtain a $2m+1$ order method.

\begin{figure}   
	\centering
	\begin{adjustbox}{max width=0.9\textwidth,center}
		\includegraphics[width=\linewidth,trim={0cm 0cm 1.75cm 0cm},clip]{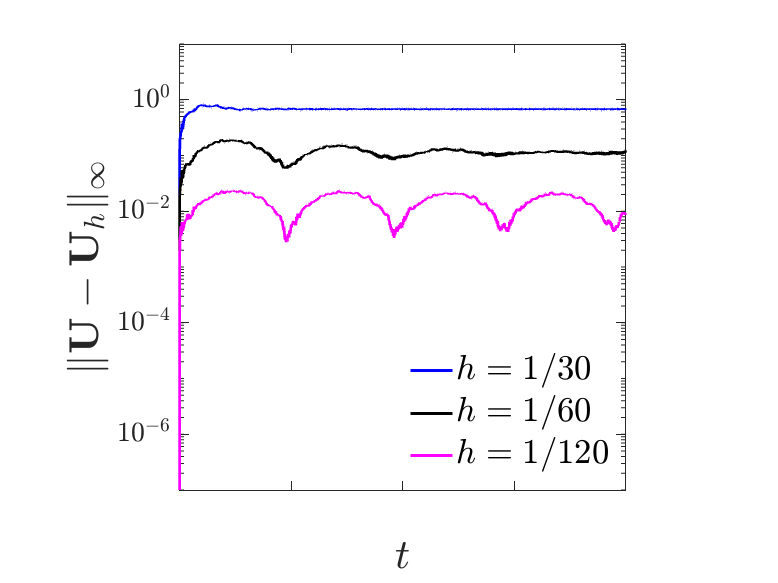} \hspace{-18.0pt}
		\includegraphics[width=\linewidth,trim={0cm 0cm 1.75cm 0cm},clip]{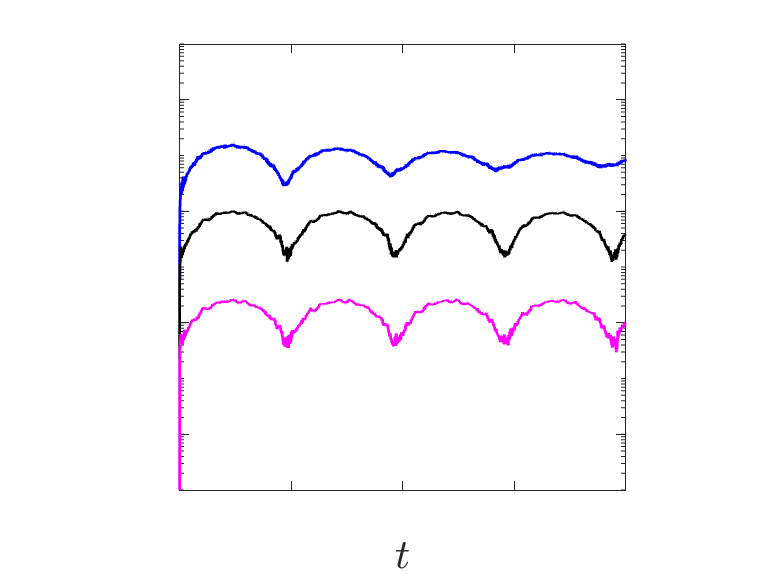}\hspace{-18pt}
	\end{adjustbox} 
       \caption{Maximum norm of the error as a function of time for a problem using generalized sheet transition  conditions \eqref{eq:GSTC_flat}, different mesh sizes and values of $m$. The left and right plots are respectively for $m=1$ and $m=2$. 
       Here $\mathbold{U} = [H_x;H_y;E_z]$.}
\label{fig:long_time_simulations_gstc_cdn}
\end{figure}

\begin{figure}   
	\centering
	\begin{adjustbox}{max width=0.9\textwidth,center}
		\includegraphics[width=\linewidth,trim={0cm 0cm 1.75cm 0cm},clip]{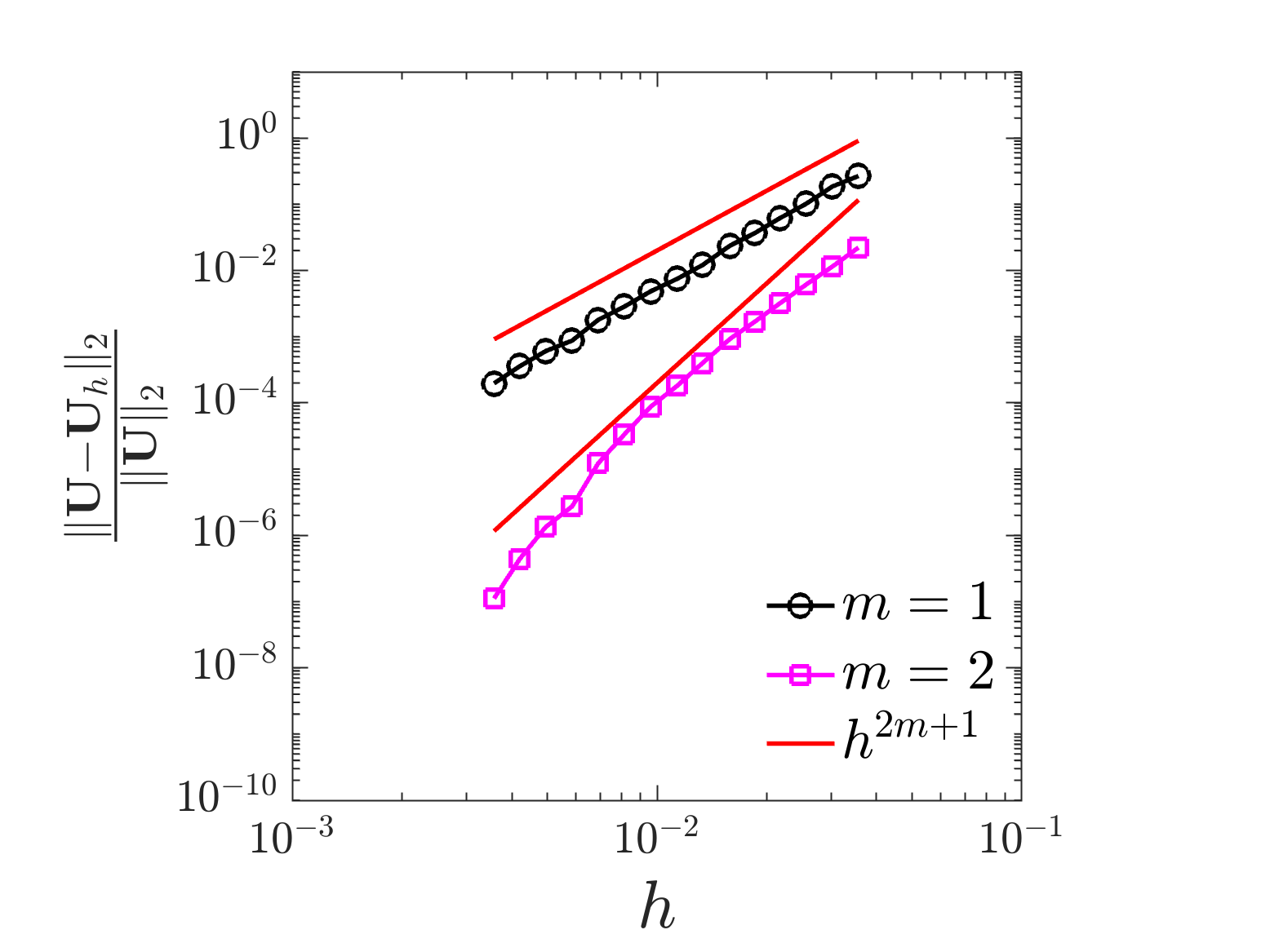} \hspace{-18.0pt}
		\includegraphics[width=\linewidth,trim={0cm 0cm 1.75cm 0cm},clip]{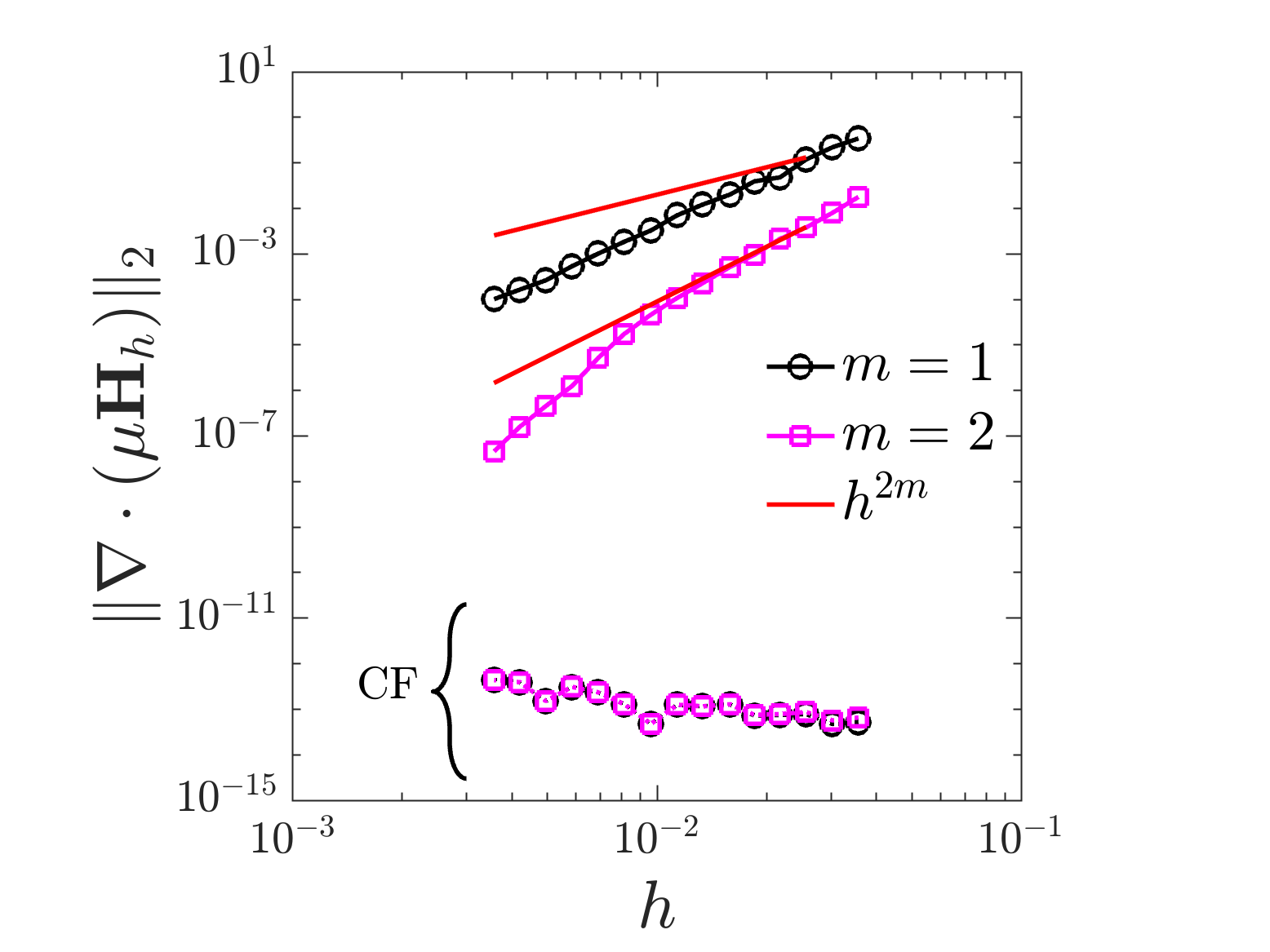}\hspace{-18pt}
	\end{adjustbox} 
       \caption{Convergence plots for a problem using GSTC model \eqref{eq:GSTC_flat} and different values of $m$. 
The left plots illustrate the relative error of $\mathbold{U} = [H_x; H_y; E_z]$ in 2-norm, 
    while the right plots show the norm \eqref{eq:computation_divergence_H} for the divergence of the magnetic field, 
    as well as the contribution of correction functions.}
\label{fig:conv_gstc}
\end{figure}

Fig.~\ref{fig:long_time_simulations_gstc_cdn} illustrates the maximum norm error as a function of time over the time interval $I=[0,100]$ with a CFL constant of $0.5$ and $m=1-2$. 
Long-time simulations suggest that the method is stable. 
Note that the Hermite-Taylor discrete correction function method becomes unstable for $m=3$. 
As before,
    we did not investigate CFL constants lower than $0.1$.
The maximum condition number for the matrices $M$ for the DCFM is around $10^5$ and $10^8$ with respectively $m=1-2$ for the PEC part, 
    while it is approximately $10^6$ and $10^8$ for the GSTC part.
Convergence plots at $t_f =1$ are shown in Fig.~\ref{fig:conv_gstc}. 
The results are in agreement with the theory and the convergence of the divergence constraint on the magnetic field is even better than expected.
To the best of our knowledge, 
    numerical methods for handling GSTC models achieve at most second-order rate of convergence \cite{Hosseini2018,Smy2020,Tian2022,Li2025}. 
The Hermite-Taylor discrete correction function method exhibits clear third and fifth-order rates of convergence and is a promising avenue for enforcing GSTC-type models. 
Future work will further investigate this direction.

\section{Conclusion}

In this work, 
    we propose a novel Hermite-Taylor correction function method 
    to enforce various conditions on a surface in an immersed setting. 
The key idea of the discrete correction function method is to directly constrain the polynomial coefficients of the correction functions approximating the electromagnetic fields when enforcing the governing equations, 
    the surface conditions 
    and the correction functions to match the Hermite-Taylor numerical solution.
This approach is devoid of any volumetric or surface quadrature rules, 
    making it more efficient than other correction function methods combined with high-order base methods. 
Moreover, 
    the resulting minimization problem is solved using a least squares approach with a pivoted QR factorization, 
    which considerably reduces the condition number of the DCFM matrices when compared to the CFM proposed in \cite{LawAppeloHagstrom2025}. 
For boundary and interface conditions, 
    the proposed method achieves up to seventh-order convergence ($m=1-3$) with a reasonable choice of the CFL constant,
    even for problems with variable coefficients and discontinuous solutions at the interface. 
For generalized sheet transition conditions, 
    we achieve up to a fifth-order convergence ($m=1-2$).
Long-time simulations suggest that this high-order method is stable for all considered cases. 
Future work will focus on extending the discrete correction function method to other types of numerical methods, 
    as well as to enforcing more complex mathematical models on surfaces.

\section*{Funding} 
We acknowledge the support of the Natural Sciences and Engineering Research Council of Canada (NSERC), [Grant No. RGPIN-2025-05071], 
    as well as the support from Fonds de recherche du Qu\'ebec (FRQ), [Grant No. 360514].
	
\section*{Declarations}

\section*{Conflicts of interest/Competing interests}
The author states that there is no conflict of interest.


\bibliographystyle{elsarticle-num}
\bibliography{meing}

\begin{thebibliography}{10}
\expandafter\ifx\csname url\endcsname\relax
  \def\url#1{\texttt{#1}}\fi
\expandafter\ifx\csname urlprefix\endcsname\relax\def\urlprefix{URL }\fi
\expandafter\ifx\csname href\endcsname\relax
  \def\href#1#2{#2} \def\path#1{#1}\fi

\bibitem{Yee1966}
K.~S. Yee, Numerical solution of initial boundary value problems involving
  {M}axwell's equations in isotropic media, {IEEE} {T}rans. {A}ntennas
  {P}ropag. 14~(3) (1966) 302--307.

\bibitem{Xie2002}
Z.~Xie, C.-H. Chan, B.~Zhang, An explicit fourth-order staggered
  finite-difference time-domain method for {M}axwell's equations, {J}.
  {C}omput. {A}ppl. {M}ath. 147 (2002) 75--98.

\bibitem{Cockburn2001}
B.~Cockburn, C.-W. Shu, Runge--kutta discontinuous {G}alerkin methods for
  convection-dominated problems, J. {S}ci. {C}omput. 16~(3) (2001) 173--261.

\bibitem{Hesthaven2002}
J.~S. Hesthaven, T.~Warburton, Nodal high-order methods on unstructured grids:
  {I}. time-domain solution of {M}axwell's equations, J. {C}omput. {P}hys. 181
  (2002) 186--221.

\bibitem{Yang1997}
B.~Yang, D.~Gottlieb, J.~S. Hesthaven, Spectral simulations of electromagnetic
  wave scattering, J. {C}omput. {P}hys. 134 (1997) 216--230.

\bibitem{Fan2002}
G.-X. Fan, Q.~H. Liu, J.~S. Hesthaven, Multidomain pseudospectral time-domain
  simulations of scattering by objects buried in lossy media, {IEEE} {T}rans.
  {G}eosci. {R}emote {S}ens. 40 (2002) 1366--1373.

\bibitem{Goodrich2005}
J.~Goodrich, T.~Hagstrom, J.~Lorenz, Hermite methods for hyperbolic
  initial-boundary value problems, Math. {C}omp. 75 (2005) 595--630.

\bibitem{Kreiss1972}
H.-O. Kreiss, J.~Oliger, Comparison of accurate methods for the integration of
  hyperbolic equations, Tellus 24~(3) (1972) 199--215.

\bibitem{Henshaw1994}
W.~D. Henshaw, H.-O. Kreiss, L.~G.~M. Reyna, A fourth-order-accurate difference
  approximation for the incompressible {N}avier-{S}tokes equations, {C}omput.
  {F}luids 23 (1994) 575--593.

\bibitem{Henshaw2006}
W.~D. Henshaw, A high-order accurate parallel solver for {M}axwell's equations
  on overlapping grids, {SIAM} {J}. {S}ci. {C}omput. 28~(5) (2006) 1730--1765.

\bibitem{Lyon2010}
M.~Lyon, O.~P. Bruno, High-order unconditionally stable {FC-AD} solvers for
  general smooth domain {II. E}lliptic, parabolic and hyperbolic {PDE}s;
  theoretical considerations, J. {C}omput. {P}hys. 229 (2010) 3358--3381.

\bibitem{Zhao2004}
S.~Zhao, G.~W. Wei, High-order {FDTD} methods via derivative matching for
  {M}axwell's equations with material interfaces, J. {C}omput. {P}hys. 200
  (2004) 60--103.

\bibitem{Zhang2016}
Y.~Zhang, D.~D. Nguyen, K.~Du, J.~Xu, S.~Zhao, Time-domain numerical solutions
  of {M}axwell interface problems with discontinuous electromagnetic waves,
  Adv. {A}ppl. {M}ath. {M}ech. 8 (2016) 353--385.

\bibitem{Dye2014}
E.~Dye, Performance analysis and optimization of {H}ermite methods on {NVIDIA
  GPU}s using {CUDA}, Master's thesis, University of {N}ew {M}exico (2014).

\bibitem{Vargas2017}
A.~Vargas, J.~Chan, T.~Hagstrom, T.~Warburton, {GPU} acceleration of {H}ermite
  methods for the simulation of wave propagation, in: Spectral and High Order
  Methods for Partial Differential Equations {ICOSAHOM} 2016, Springer
  {I}nternational {P}ublishing, Cham, 2017, pp. 357--368.

\bibitem{Appelo2018}
D.~Appel\"{o}, T.~Hagstrom, A.~Vargas, Hermite methods for the scalar wave
  equation, {SIAM} {J}. {S}ci. {C}omput. 40 (2018) A3902--A3927.

\bibitem{Vargas2019}
A.~Vargas, T.~Hagstrom, J.~Chan, T.~Warburton, Leapfrog time-stepping for
  {H}ermite methods, J. {S}ci. {C}omput. 80 (2019) 289--314.

\bibitem{Chen2012}
R.~Chen, T.~Hagstrom, {$P$}-adaptive {H}ermite methods for initial value
  problems, {ESAIM}: {M}2{AN} 46 (2012) 545--557.

\bibitem{AppeloHagstromLaw2025}
D.~Appel\"{o}, T.~Hagstrom, Y.-M. Law, Energy-conserving {H}ermite methods for
  {M}axwell's equations, Commun. {A}ppl. {M}ath. {C}omput. 7 (2025) 1146--1173.

\bibitem{Hagstrom2014}
T.~Hagstrom, D.~Appel\"{o}, Solving {PDE}s with {H}ermite interpolation, in:
  Spectral and {H}igh {O}rder {M}ethods for {P}artial {D}ifferential
  {E}quations {ICOSAHOM} 2014, 2014, pp. 31--49.

\bibitem{Chen2014}
X.~Chen, D.~Appel\"{o}, T.~Hagstrom, A hybrid {H}ermite-discontinuous
  {G}alerkin method for hyperbolic systems with application to {M}axwell's
  equations, J. {C}omput. {P}hys. 257 (2014) 501--520.

\bibitem{Beznosov2021}
O.~Beznosov, D.~Appel\"{o}, Hermite - discontinuous {G}alerkin overset grid
  methods for the scalar wave equation, Commun. {A}ppl. {M}ath. {C}omput. 3
  (2021) 391--418.

\bibitem{Loya2025}
A.~A. Loya, D.~Appel\"{o}, W.~D. Henshaw, {H}igh order accurate {H}ermite
  schemes on curvilinear grids with compatibility boundary conditions, J.
  {C}omput. {P}hys. 522 (2025) 113597.

\bibitem{LawAppelo2025}
Y.-M. Law, D.~Appel\"{o}, The {H}ermite-{T}aylor correction function method for
  {M}axwell's equations, Commun. {A}ppl. {M}ath. {C}omput. 7 (2025) 347--371.

\bibitem{LawAppeloHagstrom2025}
Y.-M. Law, D.~Appel\"{o}, T.~Hagstrom, The {H}ermite-{T}aylor correction
  function method for embedded boundary and {M}axwell's interface problems, J.
  {C}omput. {P}hys. 537 (2025) 114111.

\bibitem{Marques2011}
A.~N. Marques, J.-C. Nave, R.~R. Rosales, A correction function method for
  {P}oisson problems with interface jump conditions, J. {C}omput. {P}hys. 230
  (2011) 7567--7597.

\bibitem{Marques2017}
A.~N. Marques, J.-C. Nave, R.~R. Rosales, High order solution of {P}oisson
  problems with piecewise constant coefficients and interface jumps, J.
  {C}omput. {P}hys. 335 (2017) 497--515.

\bibitem{Marques2019}
A.~N. Marques, J.-C. Nave, R.~R. Rosales, Imposing jump conditions on
  nonconforming interfaces for the correction function method: a least squares
  approach, J. {C}omput. {P}hys. 397 (2019) 108869.

\bibitem{Abraham2018}
D.~S. Abraham, A.~N. Marques, J.-C. Nave, A correction function method for the
  wave equation with interface jump conditions, J. {C}omput. {P}hys. 353 (2018)
  281--299.

\bibitem{LawMarquesNave2020}
Y.-M. Law, A.~N. Marques, J.-C. Nave, Treatment of complex interfaces for
  {M}axwell's equations with continuous coefficients using the correction
  function method, J. {S}ci. {C}omput. 82~(3) (2020) 56.

\bibitem{LawNave2021}
Y.~M. Law, J.-C. Nave, {FDTD} schemes for {M}axwell's equations with embedded
  perfect electric conductors based on the correction function method, J.
  {S}ci. {C}omput. 88~(3) (2021) 72.

\bibitem{LawNave2022}
Y.~M. Law, J.-C. Nave, High-order {FDTD} schemes for {M}axwell's interface
  problems with discontinuous coefficients and complex interfaces based on the
  correction function method, J. {S}ci. {C}omput. 91~(1) (2022) 26.

\bibitem{Zhou2024}
H.~Zhou, W.~Ying, A correction function-based kernel-free boundary integral
  method for elliptic {PDE}s with implicitly defined interfaces, J. {C}omput.
  {P}hys. 496 (2024) 112545.

\bibitem{Assous2018}
F.~Assous, P.~Ciarlet, S.~Labrunie, Mathematical foundations of computational
  electromagnetism, Springer International Publishing, 2018.

\bibitem{Lindell2019}
I.~Lindell, A.~Sihvola, Boundary conditions in electromagnetics, John Wiley \&
  Sons, 2019.

\bibitem{Achouri2021}
K.~Achouri, C.~Caloz, Electromagnetic Metasurfaces: Theory and Applications,
  John Wiley \& Sons, 2021.

\bibitem{Gabbard2023}
J.~Gabbard, W.~M. van Rees, A high-order 3{D} immersed interface finite
  difference method for the advection-diffusion equation, in: {AIAA SCITECH}
  2023 {F}orum, 2023, p. 2480.

\bibitem{Abraham2017}
D.~S. Abraham, D.~D. Giannacopoulos, A parallel implementation of the
  correction function method for {P}oisson's equation with immersed surface
  charges, {I}{E}{E}{E} {T}rans. {M}agn. 53~(6) (2017).

\bibitem{Umashankar1993}
K.~Umashankar, A.~Taflove, Computational electrodynamics, Artech House, 1993.

\bibitem{Cai2003}
W.~Cai, S.~Deng, An upwinding embedded boundary method for {M}axwell's
  equations in media with material interfaces: 2{D} case, J. {C}omput. {P}hys.
  190 (2003) 159--183.

\bibitem{Hosseini2018}
K.~Hosseini, Z.~Atlasbaf, {PLRC-FDTD} modeling of general {GSTC}-based
  dispersive bianisotropic metasurfaces, IEEE {T}rans. {A}ntennas {P}ropag.
  66~(1) (2018) 262--270.

\bibitem{Smy2020}
T.~J. Smy, S.~A. Stewart, J.~G.~N. Rahmeier, S.~Gupta, {FDTD} simulation of
  dispersive metasurfaces with {L}orentzian surface susceptibilities, {IEEE}
  {A}ccess 8 (2020) 83027--83040.

\bibitem{Tian2022}
S.~Tian, K.~Wu, Q.~Ren, Modeling of metasurfaces using discontinuous {G}alerkin
  time-domain method based on generalized sheet transition conditions, IEEE
  {T}rans. {A}ntennas {P}ropag. 70~(8) (2022) 6905--6917.

\bibitem{Li2025}
C.~Li, Y.~Huang, J.~Li, Developing and analyzing a {DG} method for modeling of
  metasurfaces, J. {C}omput. {P}hys. 534 (2025) 114011.

\end{thebibliography}







\end{document}